\documentclass[11pt,reqno]{amsart}
\usepackage[T1]{fontenc}
\usepackage{lmodern}
\usepackage[a4paper,textwidth=6.5in,textheight=9in,hcentering]{geometry}
\usepackage{amsmath,amssymb,amsthm,mathtools,mathrsfs,bm}
\usepackage{booktabs,tabularx,array,enumitem,microtype}
\usepackage{xcolor}
\usepackage{hyperref}

\newif\ifrevision
\revisiontrue

\hypersetup{hidelinks,
 pdftitle={Entropy and frame measures for self-similar measures},
 pdfauthor={Zhiyong Wang and Wei Tang}}
\numberwithin{equation}{section}
\newtheorem{theorem}{Theorem}[section]
\newtheorem{lemma}[theorem]{Lemma}
\newtheorem{proposition}[theorem]{Proposition}
\newtheorem{corollary}[theorem]{Corollary}
\newtheorem{conjecture}[theorem]{Conjecture}
\theoremstyle{definition}

\newtheorem{example}[theorem]{Example}
\theoremstyle{remark}
\newtheorem{remark}[theorem]{Remark}
\newcommand{\R}{\mathbb R}
\newcommand{\Z}{\mathbb Z}
\newcommand{\N}{\mathbb N}
\newcommand{\Q}{\mathbb Q}
\newcommand{\C}{\mathbb C}

\newcommand{\E}{\mathbb E}
\newcommand{\PP}{\mathbb P}
\newcommand{\ee}{\mathrm e}
\newcommand{\ind}{\mathbf 1}
\newcommand{\BD}{\mathrm{BD}}
\DeclareMathOperator{\supp}{supp}

\setlist[enumerate,1]{label=\textup{(\roman*)},leftmargin=2.1em,itemsep=3pt,topsep=5pt}
\allowdisplaybreaks[1]

\title[Entropy and frame measures]{Entropy and frame measures for self-similar measures}
\author[Z. Wang]{Zhiyong Wang}
\address{School of Mathematics and Statistics, Hunan First Normal University, Changsha, Hunan 410205, P.~R.~China}
\email{wzyzzql@163.com}
\author[W. Tang]{Wei Tang$^*$}
\address{School of Mathematics and Statistics, Hunan First Normal University, Changsha, Hunan 410205, P.~R.~China}
\email{twmath2016@163.com}
\thanks{The first author is supported by the Natural Science Foundation of Hunan Province (Grant No.~2024JJ4017) and the Science and Technology Innovation Program of Hunan Province (Grant No.~2024RC3227). The second author is supported by the National Natural Science Foundation of China (Grant No.~12271156) and the Natural Science Foundation of Hunan Province (Grant No.~2023JJ40201). Corresponding author: Wei Tang.}
\date{\today}

\begin{document}
\begin{abstract}
For $0<\rho<1$ and a set of $N\ge2$ distinct real digits
$\mathscr D=\{d_0,\ldots,d_{N-1}\}\subset\mathbb{R}$, let
\[
 \mu_{\rho,\mathscr D}=\frac1N\sum_{i=0}^{N-1}\mu_{\rho,\mathscr D}\circ S_i^{-1},
 \qquad S_i(x)=\rho(x+d_i).
\]
Write $X=\sum_{j\ge1}\rho^j d_{\omega_j}$, where the symbols $\omega_j$
are independent and uniform on $\{0,\ldots,N-1\}$. We prove that $\mu_{\rho,\mathscr D}$ admits no frame measure
if the conditional entropy $H(\omega_1\mid X)>0$. In particular, this holds for all $\rho>1/N$, without any
absolute-continuity assumption or arithmetic restriction
on the digits or the contraction. At $\rho=1/N$, for commensurable digits, a frame measure exists if and only if $\mu_{\rho,\mathscr D}$ is normalized Lebesgue measure on its attractor. For $\mathscr D=\{0,\ldots,N-1\}$ and $\rho<1/N$, we also prove nonexistence whenever $\rho^{-m}=a/b$ for some $m\in\N$, with $a/b$ in lowest terms, $b\mid N$ and $N/b\nmid a$, extending the integer-power obstruction. 
The extension compares mask quotients using finite families of target pole classes in place
of one target class. 
\end{abstract}
\maketitle
\section{Introduction}\label{sec:intro}

Let $\nu$ be a compactly supported Borel probability measure on $\R$. For $f\in L^2(\nu)$, write
\[
 \widehat{f\nu}(\lambda)=\int_{\R}f(x)\ee^{-2\pi i\lambda x}\,d\nu(x),
 \qquad e_\lambda(x)=\ee^{2\pi i\lambda x}.
\]
A positive Borel measure $\Gamma$ on $\R$ is called a \emph{frame measure} for $\nu$ if there are constants $0<A\le C<\infty$ such that
\begin{equation}\label{eq:frame}
 A\|f\|_{L^2(\nu)}^2\le
 \int_{\R}|\widehat{f\nu}(\lambda)|^2\,d\Gamma(\lambda)
 \le C\|f\|_{L^2(\nu)}^2
 \qquad(f\in L^2(\nu)).
\end{equation}
The upper inequality alone defines a \emph{Bessel measure}. When $\Gamma=\sum_{\lambda\in\Lambda}\delta_\lambda$ for a countable set $\Lambda\subset\R$, \eqref{eq:frame} is the Fourier-frame inequality for the exponentials $\{e_\lambda:\lambda\in\Lambda\}$; positive weights on the atoms give weighted Fourier frames. A probability measure is \emph{spectral} if its $L^2$-space has an orthonormal basis of exponentials; the set $\Lambda$
is called a \emph{spectrum}. Counting measure on a spectrum is then a frame measure with bounds $A=C=1$. Frame measures, introduced in \cite{DHW}, allow both atomic and non-atomic frequency measures. We refer to \cite{Christensen} for general frame theory.

The existence of a frame measure imposes a restriction on overlap
that is stronger than the absence of exact overlaps. We study this
restriction for equal-weight self-similar measures on the real line. Our main result is an obstruction expressed in terms of conditional entropy. It excludes every frame measure throughout the supercritical range $\rho>1/N$ for arbitrary sets of $N$ distinct real digits.	Fix $N\ge2$ and $0<\rho<1$, and put $D_N=\{0,1,\ldots,N-1\}$. Consider
\[
 \mathscr D=\{d_0<d_1<\cdots<d_{N-1}\}\subset\R,
\]
and denote by $\mu_{\rho,\mathscr D}$ the uniform self-similar probability measure determined by
\begin{equation}\label{eq:general-ifs}
 \mu_{\rho,\mathscr D}=\frac1N\sum_{i=0}^{N-1}(S_i)_{\#}\mu_{\rho,\mathscr D},
 \qquad S_i(x)=\rho(x+d_i),
\end{equation}
where $(S_i)_{\#}\mu(B)=\mu(S_i^{-1}(B))$ for Borel sets $B\subset\R$ defines the pushforward measure. Here and below, the symbol $i\in D_N$ labels the actual digit $d_i\in\mathscr D$. On the uniform product space $D_N^{\N}$, the corresponding random series is
\[
 X=\sum_{j\ge1}\rho^j d_{\omega_j}.
\]

For absolutely continuous measures, Lai \cite[Theorems~1.1 and 1.2]{Lai11} characterized
the existence of Fourier frames by finite upper and positive lower
bounds on the density, holding almost everywhere with respect to the measure. In particular, an absolutely continuous uniform self-similar
measure with $N$ maps can admit a Fourier frame only at the critical
contraction $\rho=1/N$. This excludes almost every supercritical
Bernoulli convolution with two digits. Dutkay and Lai \cite[Theorems~2.2, 2.4 and 4.1]{DL}
established a general comparison principle for translates of
restricted measures and extended the density restrictions to frame
measures. Their comparison principle applies to singular measures
as well, but their corresponding restriction on the contraction
assumes absolute continuity.

We remove the absolute-continuity assumption from the supercritical
nonexistence conclusion. The obstruction is expressed by the
conditional entropy $H(\omega_1\mid X)$, which measures the uncertainty
of the first symbol after the spatial point is known. The definition
is recalled in Section~\ref{sec:probability}.

\begin{theorem}\label{thm:general-digits}
Let $\mathscr D\subset\R$ consist of $N\ge2$ distinct digits, let $0<\rho<1$, and put $\mu=\mu_{\rho,\mathscr D}$. 
If $H(\omega_1\mid X)>0$, then $\mu$ admits no frame measure. Moreover, $H(\omega_1\mid X)=0$ if and only if the first-level component measures
$\mu_i=N^{-1}(S_i)_{\#}\mu$, $0\le i<N$, are mutually singular.
\end{theorem}

The theorem concerns the first-level component measures, not merely the intersection of their closed supports. It applies to arbitrary real digits, whether or not there is an exact overlap. Its principal consequence is the following supercritical nonexistence result.

\begin{theorem}\label{thm:overlap}
Let $\mathscr D\subset\R$ consist of $N\ge2$ distinct digits. If $1/N<\rho<1$, then $\mu_{\rho,\mathscr D}$ admits no frame measure. 
\end{theorem}

Thus no supercritical uniform self-similar measure on the line
admits a frame measure, whether it is absolutely continuous or singular.
There is no arithmetic assumption on $\rho$ or on $\mathscr D$,
and the attractor need not be an interval. The conclusion applies both with and without exact overlaps.

The proof of the supercritical result has two ingredients.
First, we show that positive conditional entropy of the first digit
excludes a frame measure. The argument is based on
translated restrictions of the measure and conditional prefix
probabilities. Second, we prove the entropy bound
\begin{equation}\label{eq:entropy-main}
 H(\omega_1\mid X)\ge\log(N\rho)
\end{equation}
by Gaussian smoothing and elementary information inequalities.
This gives a direct proof of Theorem~\ref{thm:overlap} without using
dimension theory. The same entropy bound also follows from the
stronger projection-entropy dimension formula of Feng and Hu
\cite{FH}; see Remark~\ref{rem:FH}.

We next apply the entropy obstruction at the critical contraction,
and then study a different obstruction in the separated
consecutive-digit case. The first argument uses the uncertainty of the first digit given the
spatial point. The second uses arithmetic relations among Fourier
zeros in a setting where the coding is injective.

At the critical contraction $\rho=1/N$, the digits still matter. We call
$\mathscr D$ \emph{commensurable} if
$\mathscr D\subset a+h\Z$ for some $a\in\R$ and $h>0$.
An affine change of digits then makes all digits integers.
An \emph{exact overlap} means that $S_u=S_v$ for two distinct
words of the same length, where
$S_w=S_{w_1}\circ\cdots\circ S_{w_n}$. Write $|K|$ for the
Lebesgue measure of a Borel set $K$. The following critical
classification is an application of the entropy obstruction.

\begin{theorem}\label{thm:critical}
Let $\mathscr D$ consist of $N\ge2$ distinct commensurable real digits,
let $\rho=1/N$, and put $\mu=\mu_{1/N,\mathscr D}$ and $K=\supp\mu$.
The following conditions are equivalent.
\begin{enumerate}
\item The measure $\mu$ admits a frame measure.
\item The iterated function system has no exact overlap.
\item For every $n\ge1$,
\begin{equation}\label{eq:critical-cardinality}
 \#\bigl(\mathscr D+N\mathscr D+\cdots+N^{n-1}\mathscr D\bigr)=N^n.
\end{equation}
\item $|K|>0$.
\item $|K|>0$ and $\mu=|K|^{-1}\ind_K\,dx$.
\item $H(\omega_1\mid X)=0$.
\end{enumerate}
Under these conditions, $d\Gamma(\lambda)=|K|\,d\lambda$ is a Parseval
frame measure for $\mu$.
\end{theorem}

An exact overlap gives positive conditional entropy. In the absence of exact overlaps, counting distinct integer prefix sums gives $|K|>0$, and Plancherel's theorem supplies a frame measure. Commensurability is used in the counting argument. The equivalence can fail without it; see
Example~\ref{ex:irrational-critical}.

We next consider consecutive digits and write
$\mu_{\rho,N}=\mu_{\rho,D_N}$, where $D_N=\{0,1,\ldots,N-1\}$. Thus
\begin{equation}\label{eq:ifs}
 \mu_{\rho,N}=\frac1N\sum_{i=0}^{N-1}(S_i)_{\#}\mu_{\rho,N},
 \qquad S_i(x)=\rho(x+i).
\end{equation}
For these measures, Dai, He and Lau \cite{DHL} obtained a complete
spectral classification. Spectrality is sufficient for a frame
measure, but nonspectrality is not an obstruction by itself:
Lai and Wang \cite{LaiWang} constructed nonspectral fractal measures
with Fourier frames.

The middle-third Cantor measure is a central example.
Jorgensen and Pedersen \cite{JP} proved that it is not spectral,
and Strichartz \cite{Strichartz} asked whether it admits a Fourier
frame. De Dios Pont, Liehr and Taylor \cite{PLT} excluded Fourier
frames for odd-base Cantor measures. Fu, Song and Wang \cite{FSW}
treated the more general relation $\rho^{-m}=B$ with $B$ an odd
integer in the two-digit case. Fu and Wang \cite{FW} proved the
corresponding obstruction for arbitrary frame measures and
consecutive $N$-digit systems: $\rho^{-m}=B\in\N$ and $N\nmid B$
exclude a frame measure when $\rho<1/N$. The precise statements
needed below are recalled in Section~\ref{sec:prelim}. The second part of the paper extends this obstruction to reduced rational powers with denominator greater than one.

\begin{theorem}\label{thm:rational}
Let $N\ge2$ and $0<\rho<1/N$. Suppose that, for some $m\in\N$,
\begin{equation}\label{eq:rational-power}
 \rho^{-m}=\frac ab,\qquad a,b\in\N,\quad \gcd(a,b)=1.
\end{equation}
If
\begin{equation}\label{eq:arithmetic-global}
 b\mid N\qquad\text{and}\qquad \frac Nb\nmid a,
\end{equation}
then $\mu_{\rho,N}$ admits no frame measure.
\end{theorem}

The case $b=1$ is exactly the theorem of Fu and Wang \cite{FW}.
If $b>1$, then $\rho^{-k}\notin\N$ for every $k\in\N$:
otherwise $(a/b)^k=(\rho^{-k})^m$ would be an integer.
Thus these parameters lie outside the integer-power result.
These new cases require composite $N$ and include $\mu_{2/9,4}$ and
$\mu_{\sqrt{2/33},4}$.
For example, when $N=4$ and $q=\rho^{-1}=9/2$, multiplication
by $q$ sends $k+1/2$ into $1/4+\Z$ when $k$ is even and into
$3/4+\Z$ when $k$ is odd. Thus one target class does not suffice.

Put $q=\rho^{-1}$. The proof retains the two-coordinate character estimates of the quotient method in \cite{FW}. We allow a finite family of target pole classes and determine
exactly when it covers all source poles after rescaling. Suppose that $d$ is the least positive integer with $q^d\in\Q$,
and write $q^d=a/b$ in lowest terms. Whenever a finite pole cover
exists, the smallest number of channels is $b$; a channel specifies
one exponent and one target pole class.
The arithmetic characterization and its equivalent pointwise
comparison are given in Theorem~\ref{thm:equivalence}; their
analytic consequence is Theorem~\ref{thm:equiv-obstruction}.
They determine precisely which parameters this finite-cover method
can exclude. Example~\ref{ex:fractional} explains the new fractional
cases, and Section~\ref{sec:remaining} records the parameters
left open.

The main results, together with the known spectral classification,
give the following consequences.

\begin{corollary}\label{cor:classification}
For $1/N\le\rho<1$, the measure $\mu_{\rho,N}$ has a frame measure if and only if $\rho=1/N$. For every integer $q\ge2$, the measure $\mu_{1/q,N}$ has a frame measure if and only if $N\mid q$; in that case it is spectral.
\end{corollary}

Section~\ref{sec:prelim} fixes the notation and recalls the probability and spectral results used below. Section~\ref{sec:overlap} proves the
entropy obstruction, gives an independent proof of the entropy bound, and
deduces the supercritical and critical conclusions.
Section~\ref{sec:rational} treats the separated
consecutive-digit case and establishes the rational-power obstruction
through finite pole covers. Section~\ref{sec:remaining} first
illustrates the scope of the entropy theorem, then applies both methods to three commensurable digits, and finally
describes the remaining consecutive-digit parameters.	

\section{Preliminaries}\label{sec:prelim}

Throughout the paper, $\N=\{1,2,\ldots\}$, all logarithms are natural,
and $0\log0=0$. We write $\nu|_E$ for the restriction of $\nu$ to a Borel set $E$,
$Z(\widehat\nu)=\{t\in\R:\widehat\nu(t)=0\}$ for the real zero set of $\widehat\nu$, and
$\dim_{\rm H}\nu=\inf\{\dim_{\rm H}E:\nu(E)=1\}$ for the Hausdorff
dimension of a probability measure, where the
infimum is over Borel sets $E$. Since $\nu$ is finite,
Cauchy-Schwarz gives $L^2(\nu)\subset L^1(\nu)$. Hence, for every $f\in L^2(\nu)$
and $\lambda\in\mathbb R$, the Fourier integral $ \widehat{f\nu}(\lambda)=\int f(x)\ee^{-2\pi i\lambda x}\,d\nu(x)$ in \eqref{eq:frame} is well defined. For a fixed Bessel measure $\Gamma$ for $\nu$, define the Fourier
energy of $f\in L^2(\nu)$ by
$$
\mathcal E_\Gamma(f)
=
\int_{\mathbb R}
|\widehat{f\nu}(\lambda)|^2\,d\Gamma(\lambda).
$$
The measure $\nu$ is omitted from the notation and will be clear from the context.

\subsection{Affine changes and translations}
The following lemma describes how Bessel and frame measures change under an affine map.
\begin{lemma}\label{lem:affine}
Let $T(x)=ux+v$, where $u\ne0$, and put $\widetilde\nu=T_{\#}\nu$ and $D_a(\lambda)=a\lambda$. Then $\Gamma$ is a Bessel measure (respectively, a frame measure) for $\nu$ if and only if $(D_{1/u})_{\#}\Gamma$ is a Bessel measure (respectively, a frame measure) for $\widetilde\nu$, with the same bounds.
\end{lemma}
\begin{proof}
Put $\widetilde\Gamma=(D_{1/u})_{\#}\Gamma$. For $f\in L^2(\widetilde\nu)$, let $g=f\circ T$. The pushforward identity gives
\[
 \|g\|_{L^2(\nu)}^{2}
 =\int_{\R}|f(Tx)|^2\,d\nu(x)
 =\int_{\R}|f(y)|^2\,d\widetilde\nu(y)
 =\|f\|_{L^2(\widetilde\nu)}^{2}.
\]
For every $\lambda\in\R$,
\[
 \widehat{f\widetilde\nu}(\lambda)
 =\int_{\R}f(Tx)\ee^{-2\pi i\lambda(ux+v)}\,d\nu(x)
 =\ee^{-2\pi iv\lambda}\widehat{g\nu}(u\lambda).
\]
Since $|\ee^{-2\pi i v\lambda}|=1$,
\begin{align*}
 \int_{\R}|\widehat{f\widetilde\nu}(\lambda)|^2\,d\widetilde\Gamma(\lambda)
 =\int_{\R}|\widehat{f\widetilde\nu}(\xi/u)|^2\,d\Gamma(\xi)
 =\int_{\R}|\widehat{g\nu}(\xi)|^2\,d\Gamma(\xi).
\end{align*}
The Bessel and frame inequalities therefore hold with the same bounds. The converse follows by applying the same argument to $T^{-1}$, since $(T^{-1})_{\#}\widetilde\nu=\nu$ and $(D_u)_{\#}\widetilde\Gamma=\Gamma$.
\end{proof}

For later use, if $\theta$ is a finite measure and $T_t(x)=x+t$, then
\begin{equation}\label{eq:translation}
 \widehat{(T_t)_{\#}\theta}(\lambda)
 =\ee^{-2\pi it\lambda}\widehat\theta(\lambda).
\end{equation}
Thus a measure and its translates have Fourier transforms with the same modulus. In Section~\ref{sec:overlap}, we apply this identity to restrictions of $\mu$ and compare the $L^2(\mu)$ norms of their densities. This comparison is used in \cite{DL} to study the uniformity of measures with Fourier frames; density bounds for absolutely continuous measures were obtained in \cite{Lai11}.

\subsection{The random-series model}
For the digit set in \eqref{eq:general-ifs}, let $\Omega=D_N^{\N}$ have its product topology and the uniform product probability measure $\PP$. 
{
The coordinate maps $\omega_j:\Omega\to D_N$ are independent and uniformly distributed on $D_N$. Let $\vartheta:\Omega\to\Omega$ be the left shift, given by $\vartheta(\omega_1,\omega_2,\ldots)=(\omega_2,\omega_3,\ldots)$. For $\omega=(\omega_1,\omega_2,\ldots)\in\Omega$ and $k\ge0$, define
}
\begin{equation}\label{eq:coding}
 \pi_{\rho,\mathscr D}(\omega)=\sum_{j\ge1}\rho^j d_{\omega_j},
 \qquad X=\pi_{\rho,\mathscr D},\qquad
 X_k=\pi_{\rho,\mathscr D}\circ\vartheta^k.
\end{equation}
We write $\pi$ when the parameters are fixed. For consecutive digits $d_i=i$, this is the usual $N$-Bernoulli coding. Since the digits are bounded and $0<\rho<1$, the series converges uniformly; its finite partial sums are continuous, and hence so is $\pi$. Separating the first $n$ terms gives
\begin{equation}\label{eq:tail-decomposition}
 X=\sum_{j=1}^n\rho^j d_{\omega_j}+\rho^n X_n,
 \qquad X=\rho(d_{\omega_1}+X_1).
\end{equation}
The random variable $X_n$ has the same law as $X$ and is independent of $(\omega_1,\ldots,\omega_n)$. Conditioning on $\omega_1$ in the second identity shows that the law of $X$ satisfies \eqref{eq:general-ifs}. Uniqueness of the self-similar probability measure \cite{Hutchinson} gives $\pi_{\#}\PP=\mu_{\rho,\mathscr D}$.

Write
\begin{equation}\label{eq:general-endpoints}
 x_- =\frac{\rho d_0}{1-\rho},\qquad
 x_+ =\frac{\rho d_{N-1}}{1-\rho},\qquad
 \Delta=d_1-d_0>0.
\end{equation}
The constant digit sequences give the extreme points $x_-$ and $x_+$.
The set $\Omega$ is compact, and every nonempty open cylinder has positive $\PP$-measure. Consequently,
\[
 K:=\supp\mu_{\rho,\mathscr D}=\pi(\Omega)\subset[x_-,x_+].
\]
Indeed, $\pi(\Omega)$ is a compact set of full measure. Conversely, if $x=\pi(\omega)$, the inverse image of any open neighborhood of $x$ contains an open cylinder and hence has positive probability. For a word $w=(w_1,\ldots,w_n)\in D_N^n$, put
\begin{equation*}
t_w=\sum_{j=1}^n\rho^j d_{w_j},\qquad
S_w=S_{w_1}\circ\cdots\circ S_{w_n},\qquad
S_w(x)=t_w+\rho^n x.
\end{equation*}
Conditioning on the first $n$ symbols gives the following decomposition of $\mu=\mu_{\rho,\mathscr D}$:
\begin{equation}\label{eq:prefix-decomp}
 \mu=\sum_{w\in D_N^n}\mu_w,\qquad
 \mu_w=N^{-n}(S_w)_{\#}\mu.
\end{equation}
Words are counted separately even when their affine maps coincide.

\subsection{Conditional probabilities and entropy}\label{sec:probability}
We recall the probability facts used in Section~\ref{sec:overlap}. Let $(\Omega,\mathcal F,\PP)$ be a probability space and $\mathcal G\subset\mathcal F$ a sub-$\sigma$-algebra. For $Y\in L^1(\PP)$, the conditional expectation $\E(Y\mid\mathcal G)$ is the $\mathcal G$-measurable integrable random variable characterized, up to almost-everywhere equality, by
\[
\E[\ind_A\E(Y\mid\mathcal G)]=\E[\ind_A Y]
\qquad(A\in\mathcal G).
\]
It exists by the Radon--Nikodym theorem. Conditional probability is defined by $\PP(A\mid\mathcal G)=\E(\ind_A\mid\mathcal G)$.
For a real-valued random variable $V$, let $\sigma(V)$ denote the $\sigma$-algebra generated by $V$. If a random variable $\Xi$ takes values in $D_N$, put
\[
 u_i(V)=\PP(\Xi=i\mid V)
 =\E(\ind_{\{\Xi=i\}}\mid\sigma(V)).
\]
These conditional probabilities form a probability vector almost
surely. For a probability vector $(p_0,\ldots,p_{N-1})$, its entropy and the conditional entropy of $\Xi$ given $V$ are defined by
\[
 \mathfrak h(p_0,\ldots,p_{N-1})
 =-\sum_{i=0}^{N-1}p_i\log p_i,\qquad
 H(\Xi\mid V)=\E\,\mathfrak h((u_i(V))_i).
\]
The entropy function is continuous, concave and nonnegative on the set of probability vectors. It vanishes exactly when one entry is $1$ and all the others are $0$. Moreover,
\begin{equation}\label{eq:entropy-bounds}
 0\le H(\Xi\mid V)\le H(\Xi)\le\log N,
\end{equation}
where $H(\Xi)=\mathfrak h((\PP(\Xi=i))_i)$ is the ordinary entropy 
of $\Xi$. These definitions do not require a density for $V$; see
\cite{CoverThomas}.

We write $I(\Xi;V)=H(\Xi)-H(\Xi\mid V)$ for mutual
information. If $Q$ is a measurable function, the tower property gives
\[
 \PP(\Xi=i\mid Q(V))
 =\E\bigl(u_i(V)\mid\sigma(Q(V))\bigr).
\]
Conditional Jensen's inequality for $\mathfrak h$, followed by
expectation, therefore yields
\begin{equation}\label{eq:data-processing}
 H(\Xi\mid Q(V))\ge H(\Xi\mid V),\qquad
 I(\Xi;Q(V))\le I(\Xi;V).
\end{equation}

We also use the upward martingale convergence theorem
\cite{Williams}. If $\{\mathcal G_k\}$ is an increasing sequence of $\sigma$-algebras and
$\mathcal G=\sigma(\bigcup_k\mathcal G_k)$, then, for every
bounded random variable $Y$,
\begin{equation}\label{eq:martingale}
 \E(Y\mid\mathcal G_k)\longrightarrow\E(Y\mid\mathcal G)
 \quad\text{almost surely and in }L^1.
\end{equation}
The Bernoulli shift $\vartheta$ on $(D_N^\N,\PP)$ preserves measure
and is ergodic. We use the Birkhoff ergodic theorem \cite{Walters}
in the form
\begin{equation}\label{eq:birkhoff}
 \frac1n\sum_{j=0}^{n-1}\Psi(\vartheta^j\omega)
 \longrightarrow\E\Psi
 \quad\PP\text{-almost surely},\qquad \Psi\in L^1(\PP).
\end{equation}
The random variables to which these results are applied in Section~\ref{sec:overlap} are bounded.

\subsection{Consecutive digits and known results}
For $\mu_{\rho,N}$ in \eqref{eq:ifs}, the following spectral
classification is due to Dai, He and Lau \cite{DHL}.

\begin{theorem}\label{thm:DHL}
For $N\ge2$ and $0<\rho<1$, the measure $\mu_{\rho,N}$ is spectral if and only if $\rho=1/q$ for an integer $q$ with $N\mid q$.
\end{theorem}

The integer-power obstruction of Fu and Wang \cite{FW} is the
starting point for Section~\ref{sec:rational}.

\begin{theorem}\label{thm:FW}
Let $N\ge2$ and $0<\rho<1/N$. If $\rho^{-m}=B\in\N$ for some integer $m\ge1$ and $N\nmid B$, then $\mu_{\rho,N}$ admits no frame measure.
\end{theorem}
The existence of a spectrum when $N\mid q$ follows from compatible digit pairs; see also \cite{LW02,DHLai}. The proof of Theorem~\ref{thm:rational} includes Theorem~\ref{thm:FW} as the case $b=1$. For consecutive digits, the convex hull of the support is $[0,\ell_\rho]$, where
\begin{equation*}
\ell_\rho=\frac{(N-1)\rho}{1-\rho}.
\end{equation*}
The following lemma describes the support and coding for consecutive digits.
\begin{lemma}\label{lem:coding}
For $\mu=\mu_{\rho,N}$, let $K=\supp\mu$. If $1/N\le\rho<1$, then $K=[0,\ell_\rho]$. If $0<\rho<1/N$, then $\pi:\Omega\to K$ is a homeomorphism, and
\[
 U:L^2(\mu)\longrightarrow L^2(\PP),\qquad Uf=f\circ\pi,
\]
is unitary. Thus every $G\in L^2(\PP)$ determines a unique $f\in L^2(\mu)$ with $G=f\circ\pi$ almost surely.
\end{lemma}
\begin{proof}
If $\rho\ge1/N$, then $\ell_\rho\ge1$, so the intervals
\[
 S_i([0,\ell_\rho])=[\rho i,\rho i+\rho \ell_\rho],\qquad 0\le i<N,
\]
meet or overlap in order. The left endpoint of the first interval is $0$, and the right endpoint of the last is $\rho(N-1)+\rho \ell_\rho=\ell_\rho$. Hence $\bigcup_{i=0}^{N-1}S_i([0,\ell_\rho])=[0,\ell_\rho]$. Uniqueness of the compact invariant set gives $K=[0,\ell_\rho]$. 

Suppose that $\rho<1/N$. Let $\omega\ne\omega'$ and
let $k$ be the first coordinate at which they differ. Since the digits are distinct integers,
\[
 |\pi(\omega)-\pi(\omega')|
 \ge\rho^k-\frac{(N-1)\rho^{k+1}}{1-\rho}
 =\rho^k\frac{1-N\rho}{1-\rho}>0.
\]
Thus $\pi$ is injective. As a continuous bijection from the compact space $\Omega$ onto the Hausdorff space $K$, it is a homeomorphism.

The pushforward identity gives
\[
 \|Uf\|_{L^2(\PP)}^2=\int_\Omega|f(\pi(\omega))|^2\,d\PP(\omega)
 =\int_K|f|^2\,d\mu=\|f\|_{L^2(\mu)}^2.
\]
This also proves that $U$ is well defined on almost-everywhere equivalence classes. Since $U$ is linear, it is an isometry.

To prove surjectivity, let $G\in L^2(\PP)$ and choose a
Borel representative. Define
$$
f(x)=
\begin{cases}
G(\pi^{-1}(x)), & x\in K,\\
0, & x\notin K.
\end{cases}
$$
Since $\pi^{-1}$ is continuous and $K$ is closed, $f$ is Borel measurable. Moreover,
$$
\int_{\mathbb R}|f(x)|^2\,d\mu(x)
=
\int_\Omega |G(\omega)|^2\,d\PP(\omega)<\infty.
$$
Hence $f\in L^2(\mu)$ and $Uf=G$ $\PP$-almost everywhere.
Thus $U$ is a surjective linear isometry and hence unitary.
\end{proof}

\section{The entropy obstruction and its consequences}\label{sec:overlap}

Throughout this section, $\mathscr D=\{d_0<\cdots<d_{N-1}\}$ is an arbitrary finite real digit set, $0<\rho<1$, and $\mu=\mu_{\rho,\mathscr D}$. In the coding \eqref{eq:coding}, the symbol $i\in D_N$ represents the digit $d_i$. We impose no separation condition. The conditional probabilities below are defined by
Radon--Nikodym derivatives and are unique up to $\mu$-null sets. With the notation in \eqref{eq:general-endpoints}, set
\begin{equation}\label{eq:endpoint-main}
 E=[x_-,x_-+\Delta),\qquad F_n=[x_-,x_-+\rho^n\Delta),\qquad
 c=\mu_{\rho,\mathscr D}(E),\qquad
 f_n=\frac{N^{n/2}}{\sqrt c}\ind_{F_n}.
\end{equation}
Lemma~\ref{lem:endpoint} below proves that $c>0$ and that the functions $f_n$ have norm one.

\subsection{Endpoint restrictions and translated measures}
The component measure $\mu_w$ in \eqref{eq:prefix-decomp} records the distribution of $X$ on the event that the first $n$ symbols form the word $w$.

\begin{lemma}\label{lem:prefix-probabilities}
For every word $w\in D_N^n$, the measure $\mu_w$ has a density $p_w=d\mu_w/d\mu$ such that
$0\le p_w\le1$ and
$\sum_{w\in D_N^n}p_w=1$ $\mu$-almost everywhere. Moreover, $p_w(X)=\PP(\omega|_n=w\mid X)
\;\PP$-almost surely,
where $\omega|_n=(\omega_1,\ldots,\omega_n)$.	
\end{lemma}

\begin{proof}
Fix $n\ge1$. By \eqref{eq:prefix-decomp},
we have $\mu=\sum_{w\in D_N^n}\mu_w,$
and hence $0\le\mu_w\le\mu$ for every $w\in D_N^n$.
The Radon--Nikodym theorem gives a density
$p_w=d\mu_w/d\mu$ with $0\le p_w\le1$
$\mu$-almost everywhere. Summing the density identities
and using uniqueness of Radon--Nikodym derivatives, we obtain
$\sum_{w\in D_N^n}p_w=1$ $\mu$-almost everywhere.
We choose Borel versions of these densities taking values
in $[0,1]$ and set them equal to zero off $K=\supp\mu$.
For the empty word, we set $p_\varnothing\equiv1$.

Fix $w\in D_N^n$. On the event $\{\omega|_n=w\}$,
the tail decomposition \eqref{eq:tail-decomposition} gives $X=S_w(X_n)$.
The prefix $\omega|_n$ is uniform on $D_N^n$ and independent
of $X_n$, whose distribution is $\mu$.
Therefore, for every Borel set $B\subset\R$,
\begin{align*}
\E\bigl[\ind_{\{\omega|_n=w\}}\ind_B(X)\bigr]
&=N^{-n}\int_{\R}\ind_B(S_w(y))\,d\mu(y)\\
&=\int_B p_w(x)\,d\mu(x)
=\E\bigl[p_w(X)\ind_B(X)\bigr].
\end{align*}
The random variable $p_w(X)$ is bounded and
$\sigma(X)$-measurable. As $B$ ranges over the Borel subsets of $\R$, the events $\{X\in B\}$
form $\sigma(X)$. The preceding identity therefore verifies the defining property
of conditional expectation. Consequently,
$
p_w(X)
=\E\bigl(\ind_{\{\omega|_n=w\}}\mid\sigma(X)\bigr)
=\PP(\omega|_n=w\mid X)$ $\PP\text{-almost surely}.$		
\end{proof}

Although the coding may fail to be injective, points sufficiently close to the left endpoint have a prescribed initial word. The size of this window is determined by the gap $\Delta=d_1-d_0$ between the two smallest digits.

\begin{lemma}\label{lem:endpoint}
The number $c=\mu(E)$ is positive. For $n\ge1$, let $C_n=S_0^n$. Then
\begin{equation}\label{eq:endpoint-measure}
 \theta_n:=\mu|_{F_n}
 =N^{-n}(C_n)_{\#}(\mu|_E),
 \qquad \mu(F_n)=N^{-n}c.
\end{equation}
In particular, $\|f_n\|_{L^2(\mu)}=1$.
\end{lemma}

\begin{proof}
The definition of $x_-$ gives
\begin{equation}\label{eq:centered-series}
X-x_-=\sum_{j\ge1}\rho^j(d_{\omega_j}-d_0)\ge0.
\end{equation}
Choose an integer $k\ge1$ such that
$\rho^k(x_+-x_-)<\Delta$.
On the event $\{\omega_1=\cdots=\omega_k=0\}$,
the tail decomposition gives
\[
0\le X-x_-=\rho^k(X_k-x_-)
\le\rho^k(x_+-x_-)<\Delta.
\]
Thus this event is contained in $\{X\in E\}$.
This event has probability $N^{-k}$, so
\[
c=\mu(E)=\PP(X\in E)\ge N^{-k}>0.
\]

Now fix $n\ge1$. Since $S_0(x_-)=x_-$, we have
$C_n(y)=S_0^n(y)=x_-+\rho^n(y-x_-)$ and
$C_n(E)=F_n.$
If $\omega_j\ne0$ for some $1\le j\le n$, then
$d_{\omega_j}-d_0\ge\Delta$. All terms in
\eqref{eq:centered-series} are nonnegative, so
$
X-x_-\ge\rho^j\Delta\ge\rho^n\Delta.
$
Such a sequence cannot satisfy $X\in F_n$.
Therefore $X\in F_n$ implies
$\omega_1=\cdots=\omega_n=0$.
On this prefix event, the tail decomposition gives
$X=C_n(X_n)$. Since $C_n$ is injective and $C_n(E)=F_n$,
we obtain the event identity
\[
\{X\in F_n\}
=
\{\omega_1=\cdots=\omega_n=0\}\cap\{X_n\in E\}.
\]

The prefix $\omega|_n$ is independent of $X_n$,
and $X_n$ has distribution $\mu$.
Consequently, for every Borel set $B\subset\R$,
\begin{align*}
\mu(B\cap F_n)
&=\PP(X\in B\cap F_n)\\
&=\E\bigl[
\ind_{\{\omega_1=\cdots=\omega_n=0\}}
\ind_E(X_n)\ind_B(C_n(X_n))
\bigr]\\
&=N^{-n}\int_E\ind_B(C_n(y))\,d\mu(y)\\
&=N^{-n}(C_n)_{\#}(\mu|_E)(B).
\end{align*}
This proves the measure identity in
\eqref{eq:endpoint-measure}.
Taking $B=\R$ gives
$
\mu(F_n)=N^{-n}\mu(E)=N^{-n}c.
$
By the definition of $f_n$, it follows that
\[
\|f_n\|_{L^2(\mu)}^2
=\frac{N^n}{c}\mu(F_n)=1.
\]
The use of half-open intervals makes these identities valid even if an endpoint has positive mass.
\end{proof}

For a word $w\in D_N^n$, set
$
z_w=t_w-t_{0^n}=\sum_{j=1}^n\rho^j(d_{w_j}-d_0).
$
Since $t_{0^n}=(1-\rho^n)x_-$, we have $S_w=T_{z_w}\circ C_n$, where $T_t(x)=x+t$. Define
\begin{equation}\label{eq:lambda-w}
 \theta_{n,w}=N^{-n}(S_w)_{\#}(\mu|_E)
 =(T_{z_w})_{\#}\theta_n.
\end{equation}
These are positive measures satisfying $\theta_{n,w}\le\mu_w\le\mu$. Their densities are
\begin{equation}\label{eq:translated-density}
 h_{n,w}(x)=p_w(x)\ind_E(S_w^{-1}x),
 \qquad 0\le h_{n,w}\le p_w\le1.
\end{equation}
Indeed, for a nonnegative Borel function $\varphi$,
\begin{align*}
\int\varphi(x)p_w(x)\ind_E(S_w^{-1}x)\,d\mu(x)
&=\int\varphi(x)\ind_E(S_w^{-1}x)\,d\mu_w(x)\\
&=N^{-n}\int_E\varphi(S_wy)\,d\mu(y)
=\int\varphi\,d\theta_{n,w}.
\end{align*}
Thus $g_{n,w}=h_{n,w}/\sqrt{\mu(F_n)}$ belongs to $L^2(\mu)$, and
\[
 f_n\mu=\frac{\theta_n}{\sqrt{\mu(F_n)}},\qquad
 g_{n,w}\mu=\frac{\theta_{n,w}}{\sqrt{\mu(F_n)}}.
\]

\begin{lemma}\label{lem:translation-energy}
For $n\geq 1$, let
\begin{equation}\label{eq:Rn}
R_n=\int\sum_{w\in D_N^n}p_w(x)^2\,d\mu(x).
\end{equation}
If $\Gamma$ is a Bessel measure for $\mu$ with upper bound $C$, then	
\begin{equation}\label{eq:energy-Rn}
\mathcal E_\Gamma(f_n)\le\frac Cc R_n.
\end{equation}
\end{lemma}

\begin{proof}
Fix $n\ge1$. For each $w\in D_N^n$, the definitions of
$f_n$ and $g_{n,w}$, together with \eqref{eq:lambda-w}, give
$
g_{n,w}\mu=(T_{z_w})_{\#}(f_n\mu).
$
By \eqref{eq:translation},
$\widehat{g_{n,w}\mu}(\lambda)
=\ee^{-2\pi iz_w\lambda}\widehat{f_n\mu}(\lambda)$
$(\lambda\in\R).$
Since the phase factor has modulus one, we have
$
\mathcal E_\Gamma(g_{n,w})=\mathcal E_\Gamma(f_n).
$
Applying the Bessel inequality to $g_{n,w}$ and using
$g_{n,w}=h_{n,w}/\sqrt{\mu(F_n)}$, we obtain
\[
\mathcal E_\Gamma(f_n)
=\mathcal E_\Gamma(g_{n,w})
\le C\|g_{n,w}\|_{L^2(\mu)}^2
=\frac{C}{\mu(F_n)}
\int h_{n,w}(x)^2\,d\mu(x).
\]

There are $N^n$ words in $D_N^n$, and the left-hand side
does not depend on $w$. Summing over these words and using
$0\le h_{n,w}\le p_w$ from \eqref{eq:translated-density},
we obtain
\begin{align*}
N^n\mathcal E_\Gamma(f_n)
&\le
\frac{C}{\mu(F_n)}
\sum_{w\in D_N^n}\int h_{n,w}(x)^2\,d\mu(x)\\
&\le
\frac{C}{\mu(F_n)}
\int\sum_{w\in D_N^n}p_w(x)^2\,d\mu(x)\\
&=
\frac{C}{\mu(F_n)}R_n.
\end{align*}
By \eqref{eq:endpoint-measure}, $\mu(F_n)=N^{-n}c$.
Consequently,
\[
N^n\mathcal E_\Gamma(f_n)
\le\frac{CN^n}{c}R_n.
\]
Dividing by $N^n$ proves \eqref{eq:energy-Rn}.
Distinct words are counted separately, even when their
translation parameters coincide.		
\end{proof}

Since there are $N^n$ words and $\mu(F_n)=N^{-n}c$, the factor $C/c$ is independent of $n$. It remains to show that $R_n\to0$ when $H>0$.

\subsection{Conditional prefix probabilities}
For a one-letter word $i\in D_N$, write $p_i$ for the density in Lemma~\ref{lem:prefix-probabilities}. Thus $p_i(X)=\PP(\omega_1=i\mid X)$ almost surely. Put
\begin{equation}\label{eq:H}
 H=H(\omega_1\mid X)
 =-\int\sum_{i=0}^{N-1}p_i(x)\log p_i(x)\,d\mu(x).
\end{equation}
By \eqref{eq:entropy-bounds}, $0\le H\le\log N$.

\begin{lemma}\label{lem:chain}
For a symbol $i\in D_N$ and a finite word $v$,
\begin{equation}\label{eq:RN-chain}
 p_{iv}(x)=p_i(x)p_v(S_i^{-1}x)
 \quad\mu\text{-almost everywhere}.
\end{equation}
Along the random coding, simultaneously for all $n\ge1$,
\begin{equation}\label{eq:actual-chain}
 p_{\omega|_n}(X)
 =\prod_{j=0}^{n-1}p_{\omega_{j+1}}(X_j)
 \quad\PP\text{-almost surely},
\end{equation}
where $X_0=X$.
\end{lemma}

\begin{proof}
Fix $i\in D_N$ and a finite word $v$.
Since $S_{iv}=S_i\circ S_v$, the definition of the word
measures gives $\mu_{iv}=N^{-1}(S_i)_{\#}\mu_v.$
Using $d\mu_v=p_v\,d\mu$ and
$\mu_i=N^{-1}(S_i)_{\#}\mu$, we obtain, for every
nonnegative Borel function $\varphi$,
\begin{align*}
\int\varphi(x)\,d\mu_{iv}(x)
&=N^{-1}\int\varphi(S_i(y))\,d\mu_v(y)\\
&=\int\varphi(x)p_v(S_i^{-1}(x))\,d\mu_i(x)\\
&=\int\varphi(x)p_i(x)p_v(S_i^{-1}(x))\,d\mu(x).
\end{align*}
The chosen Borel versions make $p_v\circ S_i^{-1}$ measurable. Thus $p_i\cdot(p_v\circ S_i^{-1})$ is a Radon--Nikodym density
of $\mu_{iv}$ with respect to $\mu$.
Uniqueness of these densities proves \eqref{eq:RN-chain}.
For the empty word, the same identity follows from
$p_\varnothing=1$.

We next prove \eqref{eq:actual-chain}.
Since there are only countably many pairs $(i,v)$,
there is a Borel set $\mathcal N\subset\R$ with
$\mu(\mathcal N)=0$ outside which \eqref{eq:RN-chain}
holds for every such pair.
Each $X_j$, $j\ge0$, has distribution $\mu$, so
\[
\PP\left(\bigcup_{j\ge0}\{X_j\in\mathcal N\}\right)
\le\sum_{j\ge0}\PP(X_j\in\mathcal N)
=\sum_{j\ge0}\mu(\mathcal N)=0.
\]
Consequently, outside one $\PP$-null set,
\eqref{eq:RN-chain} can be applied at every $X_j$.

Fix a sequence outside this null set and an integer $n\ge1$.
The tail decomposition gives
$X_j=S_{\omega_{j+1}}(X_{j+1})$ and
$S_{\omega_{j+1}}^{-1}(X_j)=X_{j+1},$ $j\ge0.$
Applying \eqref{eq:RN-chain} at $X_j$ therefore yields
\[
p_{(\omega_{j+1},\ldots,\omega_n)}(X_j)
=
p_{\omega_{j+1}}(X_j)\,
p_{(\omega_{j+2},\ldots,\omega_n)}(X_{j+1})
\qquad(0\le j<n).
\]
When $j=n-1$, the last word is empty and its density
is $p_\varnothing=1$.
Iterating these identities from $j=0$ to $j=n-1$,
and recalling that $X_0=X$, we obtain
\[
p_{\omega|_n}(X)
=
\prod_{j=0}^{n-1}p_{\omega_{j+1}}(X_j).
\]
The exceptional set is independent of $n$, so this identity
holds $\PP$-almost surely, simultaneously for all $n\ge1$.		
\end{proof}

We apply Birkhoff's ergodic theorem to a bounded function of the first-symbol conditional probability.

\begin{lemma}\label{lem:entropy-decay}
Let $H$ and $R_n$ be defined by \eqref{eq:H} and \eqref{eq:Rn}.
If $H>0$, then $R_n\to0$.
\end{lemma}
\begin{proof}
Define
\[
 \Psi(\omega)=1-p_{\omega_1}(X(\omega)),\qquad
 \kappa=\E\Psi.
\]
Since $0\le p_i\le1$, we have $0\le\Psi\le1$. The conditional
probability identity in Lemma~\ref{lem:prefix-probabilities} gives
\[
 \kappa
 =1-\sum_{i=0}^{N-1}
       \E[\ind_{\{\omega_1=i\}}p_i(X)]
 =\int\left(1-\sum_{i=0}^{N-1}p_i(x)^2\right)\,d\mu(x).
\]
For a probability vector $(p_i)_{i=0}^{N-1}$, we have
\[
1-\sum_{i=0}^{N-1}p_i^2
=\sum_{i=0}^{N-1}p_i(1-p_i)\ge0.
\]
The sum vanishes exactly when each $p_i$ is either $0$ or $1$. Since $\sum_i p_i=1$, this means that exactly one entry is $1$, which is also the condition for the entropy to vanish. Thus $\kappa=0$ implies that $(p_i(x))_i$ has zero entropy for $\mu$-almost every $x$, and hence $H=0$. Therefore $H>0$ implies $\kappa>0$.

For $0\le t\le1$, we have $t\le\exp(-(1-t))$. The chain rule
\eqref{eq:actual-chain} therefore yields
\[
 p_{\omega|_n}(X)
 \le\exp\left(-\sum_{j=0}^{n-1}
                \Psi(\vartheta^j\omega)\right)
 \quad\PP\text{-almost surely}.
\]
By the Birkhoff theorem \eqref{eq:birkhoff}, the sum in the exponent divided by $n$
converges almost surely to $\kappa>0$. Thus
$p_{\omega|_n}(X)\to0$ almost surely.

Finally, conditioning on $X$ gives
\begin{align*}
 \E[p_{\omega|_n}(X)]
 =\sum_{w\in D_N^n}
   \E[\ind_{\{\omega|_n=w\}}p_w(X)]
 =\sum_{w\in D_N^n}\E[p_w(X)^2]
 =R_n.
\end{align*}
Since $0\le p_{\omega|_n}(X)\le1$, dominated convergence gives
$R_n\to0$.
\end{proof}

\subsection{Proof of Theorem~\ref{thm:general-digits}}

We next show that $H=0$ is equivalent to mutual singularity
of the component measures $\mu_i$.
The component measures $\mu_i$ are mutually singular if there are pairwise disjoint Borel sets
$A_0,\ldots,A_{N-1}\subset\R$ such that
$\mu_i(\R\setminus A_i)=0$ for every $0\le i<N$.

\begin{lemma}\label{lem:singular-components}
The condition $H(\omega_1\mid X)=0$ holds if and only if the measures $\mu_i=N^{-1}(S_i)_{\#}\mu$, $0\le i<N$, are mutually singular.
\end{lemma}

\begin{proof}
Suppose first that $H=0$.
Since the integrand in \eqref{eq:H} is nonnegative,
\[
\mathfrak h\bigl((p_i(x))_{i=0}^{N-1}\bigr)=-\sum_{i=0}^{N-1} p_i(x)\log p_i(x)=0
\qquad\mu\text{-almost everywhere}.
\]
By Lemma~\ref{lem:prefix-probabilities}, the vector
$(p_i(x))_{i=0}^{N-1}$ is a probability vector
$\mu$-almost everywhere.
Since the Shannon entropy of a probability vector vanishes if and only if one coordinate is $1$ and all others are $0$, it follows that, for $\mu$-almost every $x$, exactly one $p_i(x)$ equals $1$. Choose a Borel set $K_0\subset K$ with $\mu(K_0)=1$
on which these properties hold, and define
$A_i=\{x\in K_0:p_i(x)=1\},$ $0\le i<N.$
Then
\[
K_0=\bigcup_{i=0}^{N-1}A_i,\qquad
A_i\cap A_j=\varnothing\quad(i\ne j),
\]
and $p_i=\ind_{A_i}$ $\mu$-almost everywhere.
Consequently,
\[
\mu_i(\R\setminus A_i)
=\int_{\R\setminus A_i}p_i(x)\,d\mu(x)
=0.
\]
Thus the measures $\mu_i$, $0\le i<N$, are mutually
singular.

Conversely, suppose that these measures are mutually
singular. Choose pairwise disjoint Borel sets
$A_0,\ldots,A_{N-1}$ such that
$\mu_i(\R\setminus A_i)=0$ $(0\le i<N).$
Fix $i\in D_N$ and a Borel set $B\subset\R$.
For $j\ne i$, disjointness gives
$B\cap A_i\subset\R\setminus A_j$, so
$\mu_j(B\cap A_i)=0$.
Using $\mu=\sum_{j=0}^{N-1}\mu_j$, we obtain
\begin{align*}
\mu(B\cap A_i)
&=\sum_{j=0}^{N-1}\mu_j(B\cap A_i)\\
&=\mu_i(B\cap A_i)\\
&=\mu_i(B).
\end{align*}
Hence $\mu_i=\mu|_{A_i}$.
By uniqueness of Radon--Nikodym derivatives,
\[
p_i=\ind_{A_i}
\qquad\mu\text{-almost everywhere}.
\]
Every term $-p_i\log p_i$ therefore vanishes
$\mu$-almost everywhere, and \eqref{eq:H} gives $H=0$.	
\end{proof}

The equivalence in Theorem~\ref{thm:general-digits} is now established.
The nonexistence assertion follows from the energy estimate
\eqref{eq:energy-Rn} and Lemma~\ref{lem:entropy-decay}.

\begin{proof}[Proof of Theorem~\ref{thm:general-digits}]
Recall from \eqref{eq:endpoint-main} that
\begin{equation*}
E=[x_-,x_-+\Delta),\qquad F_n=[x_-,x_-+\rho^n\Delta),\qquad
c=\mu_{\rho,\mathscr D}(E),\qquad
f_n=\frac{N^{n/2}}{\sqrt c}\ind_{F_n}.
\end{equation*}
Here $x_-=\rho d_0/(1-\rho)$ and $\Delta=d_1-d_0>0$.
By Lemma~\ref{lem:endpoint}, we have $c>0$ and
$\mu(F_n)=N^{-n}c$. Hence
\[
\|f_n\|_{L^2(\mu)}^2
=\frac{N^n}{c}\mu(F_n)=1
\qquad(n\ge1).
\]

Suppose that $H(\omega_1\mid X)>0$. For every Bessel measure
$\Gamma$ with upper bound $C$, Lemmas~\ref{lem:translation-energy}
and \ref{lem:entropy-decay} imply
\[
0\le\mathcal E_\Gamma(f_n)
=\int_{\R}|\widehat{f_n\mu}(\lambda)|^2\,d\Gamma(\lambda)
\le\frac Cc R_n\longrightarrow0.
\]
If $\Gamma$ were a frame measure with lower bound $A>0$, then
\eqref{eq:frame} would give 
\[
\mathcal E_\Gamma(f_n)
\ge A\|f_n\|_{L^2(\mu)}^2
=A
\]
for every $n$, a contradiction. Thus $\mu$ admits no frame measure.
The same sequence $\{f_n\}$ satisfies $\mathcal E_\Gamma(f_n)\to0$
for every Bessel measure $\Gamma$ for $\mu$, since it depends only on $\rho$ and $\mathscr D$.

The remaining equivalence was proved in Lemma~\ref{lem:singular-components}.
\end{proof}

\subsection{An entropy lower bound}
We now prove \eqref{eq:entropy-main} directly from self-similarity. We first bound the information about $\omega_1$ contained in a Gaussian-smoothed version of $X$ and then pass to $X$ using finite partitions. The proof requires neither absolute continuity of $\mu$ nor a formula for its dimension.

Enlarging the probability space if necessary, let $Z$ be a standard Gaussian random variable independent of the digit sequence $(\omega_j)_{j\ge1}$, and set $V_\varepsilon=X+\varepsilon Z$ for $\varepsilon>0$. We keep the notation $\PP$ and $\E$ on the enlarged space. The density of $V_\varepsilon$ is
\[
 \varphi_\varepsilon(t)
 =\int_{\R}\frac1{\sqrt{2\pi}\varepsilon}
       \exp\left(-\frac{(t-x)^2}{2\varepsilon^2}\right)\,d\mu(x).
\]
We first check that its differential entropy is finite. Put
$R=\max\{|x_-|,|x_+|\}$. Since $|X|\le R$, we have
\[
 \frac{\ee^{-(|t|+R)^2/(2\varepsilon^2)}}
      {\sqrt{2\pi}\varepsilon}
 \le\varphi_\varepsilon(t)
 \le\frac1{\sqrt{2\pi}\varepsilon}.
\]
Consequently,
$|\log\varphi_\varepsilon(t)|\le C_\varepsilon(1+t^2)$, where $C_\varepsilon<\infty$ depends on $\varepsilon$
and the fixed support bound $R$.
Since $\E[V_\varepsilon^2]=\E[X^2]+\varepsilon^2\le R^2+\varepsilon^2$, it follows that
\begin{align*}
\int_{\R}\varphi_\varepsilon(t)
|\log \varphi_\varepsilon(t)|\,dt
&\le
C_\varepsilon
\int_{\R}(1+t^2)\varphi_\varepsilon(t)\,dt\\
&=C_\varepsilon\bigl(1+\E[V_\varepsilon^2]\bigr)\\
&\le C_\varepsilon(1+R^2+\varepsilon^2)
<\infty.
\end{align*}
In particular, the differential entropy
\[
 h_\varepsilon:=h(V_\varepsilon)
 =-\int_{\R}\varphi_\varepsilon(t)
                    \log\varphi_\varepsilon(t)\,dt
\]
is a finite real number. The same estimates apply to the conditional density of $V_\varepsilon$ given $\omega_1=i$.
Indeed, the conditional law of $X$ is supported in $[x_-,x_+]$, and $Z$ remains a standard Gaussian variable independent of $X$ under this conditioning.

\begin{lemma}\label{lem:Gaussian}
For every $\varepsilon>0$,
\begin{equation}\label{eq:Gaussian-information}
 I(\omega_1;V_\varepsilon)\le-\log\rho.
\end{equation}
\end{lemma}
\begin{proof}
By \eqref{eq:tail-decomposition},
$X=\rho(d_{\omega_1}+X_1)$, where $X_1$ is independent of
$\omega_1$ and has the same law as $X$. Since $Z$ is independent of the digit sequence,
the conditional law of $V_\varepsilon$ given $\omega_1=i$
is the law of
$
\rho d_i+\rho X_1+\varepsilon Z
=\rho\bigl(d_i+X_1+(\varepsilon/\rho)Z\bigr).
$
The variable $X_1+(\varepsilon/\rho)Z$ has density
$\varphi_{\varepsilon/\rho}$. Hence the conditional density
of $V_\varepsilon$ given $\omega_1=i$ is
\[
 \varphi_{\varepsilon,i}(t)
 =\rho^{-1}\varphi_{\varepsilon/\rho}(t/\rho-d_i).
\]
Changing variables $t=\rho(y+d_i)$ gives
\begin{align}\label{eq:h(V|i)}
h(V_\varepsilon\mid\omega_1=i)
&=h\left(\rho\bigl(d_i+X_1+(\varepsilon/\rho)Z\bigr)\right)\nonumber\\
&=-\int_{\R}\varphi_{\varepsilon,i}(t)
\log\varphi_{\varepsilon,i}(t)\,dt\nonumber\\
&=-\int_{\R}\varphi_{\varepsilon/\rho}(y)
\log\left(\frac{\varphi_{\varepsilon/\rho}(y)}{\rho}\right)
\,dy\nonumber\\
&=\log\rho+h_{\varepsilon/\rho}.
\end{align}
This entropy is independent of $i$. The law of total probability and Bayes' formula give
\[
 \varphi_\varepsilon=\frac1N\sum_{i=0}^{N-1}\varphi_{\varepsilon,i},
 \qquad
 \PP(\omega_1=i\mid V_\varepsilon=t)
 =\frac{\varphi_{\varepsilon,i}(t)}{N\varphi_\varepsilon(t)}.
\]
Using these identities and \eqref{eq:h(V|i)} in the definition of mutual information, we obtain
\begin{align}
I(\omega_1;V_\varepsilon)
&=\frac1N\sum_{i=0}^{N-1}\int_{\R}
\varphi_{\varepsilon,i}(t)
\log\frac{\varphi_{\varepsilon,i}(t)}
{\varphi_\varepsilon(t)}\,dt\nonumber\\
&=h_\varepsilon-\frac1N\sum_{i=0}^{N-1}
h(V_\varepsilon\mid\omega_1=i)\nonumber\\
&=h_\varepsilon-h_{\varepsilon/\rho}-\log\rho.
\label{eq:entropy-difference}
\end{align}

Let $\gamma_\varepsilon$ be the centered Gaussian density with
variance $\varepsilon^2(\rho^{-2}-1)$. Addition of independent
Gaussian variables gives
\[
 \varphi_{\varepsilon/\rho}
 =\varphi_\varepsilon*\gamma_\varepsilon.
\]

Since $u\log u$ is convex and
$\int_{\R}\gamma_\varepsilon(y)\,dy=1$,
Jensen's inequality gives
\[
\begin{aligned}
\varphi_{\varepsilon/\rho}(t)
\log\varphi_{\varepsilon/\rho}(t)
\le
\int_{\R}\varphi_\varepsilon(t-y)
\log\varphi_\varepsilon(t-y)
\gamma_\varepsilon(y)\,dy.
\end{aligned}
\]
The right-hand side is absolutely integrable in $t$, since
\[
\begin{aligned}
\int_{\R}\int_{\R}
\varphi_\varepsilon(t-y)
|\log\varphi_\varepsilon(t-y)|
\gamma_\varepsilon(y)\,dy\,dt
=
\int_{\R}\varphi_\varepsilon(t)
|\log\varphi_\varepsilon(t)|\,dt<\infty.
\end{aligned}
\]
Thus Fubini's theorem yields
\[
\int_{\R}\varphi_{\varepsilon/\rho}
\log\varphi_{\varepsilon/\rho}
\le
\int_{\R}\varphi_\varepsilon
\log\varphi_\varepsilon.
\]
Hence $h_{\varepsilon/\rho}\ge h_\varepsilon$, and
\eqref{eq:entropy-difference} proves
\eqref{eq:Gaussian-information}.
\end{proof}

We now pass to $X$ using finite partitions and their discrete entropies.

\begin{lemma}\label{lem:entropy-bound}
For any set of distinct real digits $\mathscr D$ in \eqref{eq:coding} and any $0<\rho<1$,
\[
 H(\omega_1\mid X)\ge\log(N\rho).
\]
\end{lemma}

\begin{proof}
Since $\mu$ has at most countably many atoms, choose a countable
dense set of distinct points $\{a_j:j\ge1\}\subset\R$ such that
$\mu(\{a_j\})=0$ for every $j$. For each $k\ge1$, define
\[
Q_k:\R\longrightarrow\{0,1,\ldots,k\},
\qquad
Q_k(x)=\sum_{j=1}^k\ind_{(a_j,\infty)}(x).
\]
Let $b_{k,1}<\cdots<b_{k,k}$ be the increasing rearrangement
of $a_1,\ldots,a_k$. Then
\[
Q_k^{-1}(\{r\})=
\begin{cases}
(-\infty,b_{k,1}], & r=0,\\
(b_{k,r},b_{k,r+1}], & 1\le r<k,\\
(b_{k,k},\infty), & r=k.
\end{cases}
\]
Thus $Q_k$ records the interval containing $x$, and adding
the point $a_{k+1}$ refines the preceding partition.

Fix $k$ and use the same pair $(X,Z)$ to define
$V_\varepsilon=X+\varepsilon Z$ for all $\varepsilon>0$.
Then $V_\varepsilon\to X$ almost surely as
$\varepsilon\downarrow0$. Since
$\PP(X\in\{a_1,\ldots,a_k\})=\sum_{j=1}^k\mu(\{a_j\})=0$, we have
\[
\delta_k(\omega)
:=\frac12\min_{1\le j\le k}|X(\omega)-a_j|>0
\qquad\PP\text{-almost surely}.
\]
For each such $\omega$, all sufficiently small $\varepsilon$
satisfy
$|V_\varepsilon(\omega)-X(\omega)|<\delta_k(\omega)$.
Hence $V_\varepsilon(\omega)$ and $X(\omega)$ lie in the same
partition interval, so
\[
Q_k(V_\varepsilon)=Q_k(X)
\quad\text{for all sufficiently small }\varepsilon,
\quad\PP\text{-almost surely}.
\]

For every $i\in D_N$ and $r\in\{0,\ldots,k\}$, the indicators
of the events
$\{\omega_1=i,\ Q_k(V_\varepsilon)=r\}$ therefore converge
almost surely to the indicator of
$\{\omega_1=i,\ Q_k(X)=r\}$. Since these indicators are bounded
by $1$, dominated convergence gives
\[
\PP(\omega_1=i,\ Q_k(V_\varepsilon)=r)
\longrightarrow
\PP(\omega_1=i,\ Q_k(X)=r).
\]
Both pairs take values in the same finite set
$D_N\times\{0,\ldots,k\}$. Since entropy is continuous for probability distributions on a fixed finite set, with $0\log0=0$,
\[
\begin{aligned}
H(\omega_1\mid Q_k(V_\varepsilon))
&=H(\omega_1,Q_k(V_\varepsilon))
-H(Q_k(V_\varepsilon))\\
&\longrightarrow
H(\omega_1,Q_k(X))-H(Q_k(X))\\
&=H(\omega_1\mid Q_k(X)).
\end{aligned}
\]
On the other hand, \eqref{eq:data-processing} and
Lemma~\ref{lem:Gaussian} imply
\[
\begin{aligned}
H(\omega_1\mid Q_k(V_\varepsilon))
\ge H(\omega_1\mid V_\varepsilon)
=\log N-I(\omega_1;V_\varepsilon)
\ge\log(N\rho).
\end{aligned}
\]
Letting $\varepsilon\downarrow0$, we obtain
\begin{equation}\label{eq:partition-information}
H(\omega_1\mid Q_k(X))\ge\log(N\rho)
\qquad(k\ge1).
\end{equation}

We now let $k\to\infty$. The definition of $Q_k$ gives the
$\sigma$-algebras
\[
\mathcal G_k:=\sigma(Q_k(X))
=\sigma(\{X\le a_j\}:1\le j\le k).
\]
Indeed, $Q_k(X)$ is measurable with respect to the
$\sigma$-algebra on the right, while
$
\{X\le b_{k,r}\}
=\{Q_k(X)\le r-1\}\;(1\le r\le k)
$
gives the reverse inclusion. Consequently,
$\mathcal G_k\subset\mathcal G_{k+1}$.
Moreover, density of $\{a_j:j\ge1\}$ implies that
$
\{X<t\}
=\bigcup_{\{j:a_j<t\}}\{X\le a_j\}
$ for every
$t\in\R$.
Since open half-lines generate the Borel $\sigma$-algebra
of $\R$, it follows that
\[
\sigma\left(\bigcup_{k\ge1}\mathcal G_k\right)=\sigma(X).
\]

For each $i\in D_N$, the upward martingale convergence theorem
\eqref{eq:martingale} now yields
\[
u_{k,i}:=\E(\ind_{\{\omega_1=i\}}\mid\mathcal G_k)
\longrightarrow
\E(\ind_{\{\omega_1=i\}}\mid\sigma(X))
=p_i(X)
\quad\PP\text{-almost surely}.
\]
Almost surely, $(u_{k,i})_{i=0}^{N-1}$ is a probability vector for every $k$.
Since $\mathfrak h$ is continuous on the probability simplex
and satisfies $0\le\mathfrak h\le\log N$, dominated convergence
gives
\[
\begin{aligned}
H(\omega_1\mid Q_k(X))
=\E\,\mathfrak h((u_{k,i})_{i=0}^{N-1})
\longrightarrow
\E\,\mathfrak h((p_i(X))_{i=0}^{N-1})
=H(\omega_1\mid X).
\end{aligned}
\]
Passing to the limit in \eqref{eq:partition-information},
we conclude that
$
H(\omega_1\mid X)\ge\log(N\rho).
$
\end{proof}

\begin{remark}\label{rem:FH}
The preceding entropy bound also follows from the stronger dimension
formula of Feng and Hu \cite[Definition~2.1 and Theorem~2.8]{FH}.
For the maps $S_i(x)=\rho(x+d_i)$ and the uniform Bernoulli measure, their projection entropy is $\log N-H(\omega_1\mid X)$, since $\omega_1$ is independent of $X_1$. The Lyapunov exponent is $-\log\rho$, so their formula gives
\[
 \dim_{\rm H}\mu_{\rho,\mathscr D}
 =\frac{\log N-H(\omega_1\mid X)}{-\log\rho},
\]
without separation or absolute-continuity assumptions. Since
$\dim_{\rm H}\mu_{\rho,\mathscr D}\le1$, the asserted bound follows.
The argument above provides an independent proof using Gaussian
smoothing and finite partitions, without appealing to the
dimension formula.
\end{remark}

\begin{proof}[Proof of Theorem~\ref{thm:overlap}]
For $\rho>1/N$, Lemma~\ref{lem:entropy-bound} gives
\[
 H(\omega_1\mid X)\ge\log(N\rho)>0.
\]
The conclusion follows from Theorem~\ref{thm:general-digits}.
\end{proof}

\subsection{Exact overlaps and the critical contraction}
An exact overlap gives positive conditional entropy. We shall use this fact in the proof of Theorem~\ref{thm:critical}.

\begin{proposition}\label{prop:exact}
Let $\mu=\mu_{\rho,\mathscr D}$ be as in \eqref{eq:general-ifs}.
If $S_u=S_v$ for two distinct words of the same length, then
$H(\omega_1\mid X)>0$. In particular, $\mu$ admits no frame measure.
\end{proposition}
\begin{proof}
Let $w_0$ be the longest common prefix of $u$ and $v$, and
write $u=w_0u'$ and $v=w_0v'$. We use the convention
$S_{\varnothing}=\mathrm{id}$. Then
$
S_{w_0}\circ S_{u'}
=S_u=S_v
=S_{w_0}\circ S_{v'}.
$
Since $S_{w_0}$ is injective, we obtain $S_{u'}=S_{v'}$.
The words $u'$ and $v'$ have the same length $m\ge1$.
By the choice of the longest common prefix, their first
symbols are distinct; write $u'_1=i$ and $v'_1=j$, with
$i\ne j$.

Define the positive measure
$
\eta=N^{-m}(S_{u'})_{\#}\mu
=N^{-m}(S_{v'})_{\#}\mu.
$
Since $\mu$ is a probability measure,
$\eta(\R)=N^{-m}>0$.
Iterating the self-similarity identity
\eqref{eq:general-ifs} gives, for each $\ell\in D_N$,
\[
\mu_\ell
=\frac1N(S_\ell)_{\#}\mu
=N^{-m}
\sum_{\substack{w\in D_N^m\\ w_1=\ell}}
(S_w)_{\#}\mu.
\]
All summands are positive measures. Since $u'_1=i$ and
$v'_1=j$, the measure $\eta$ occurs as a summand in both
decompositions. Hence $\eta\le\mu_i$ and $\eta\le\mu_j.$
Suppose that $\mu_i$ and $\mu_j$ were mutually singular.
Then there would exist a Borel set $E\subset\R$ such that
$\mu_i(\R\setminus E)=0$ and $\mu_j(E)=0$.
The preceding inequalities would imply
\[
0<\eta(\R)
=\eta(E)+\eta(\R\setminus E)
\le\mu_j(E)+\mu_i(\R\setminus E)
=0,
\]
a contradiction. Thus $\mu_i$ and $\mu_j$ are not mutually
singular. Lemma~\ref{lem:singular-components} yields
$H(\omega_1\mid X)>0$, and
Theorem~\ref{thm:general-digits} shows that $\mu$ admits
no frame measure.	
\end{proof}

We now prove Theorem~\ref{thm:critical}. The main step is to show that, at $\rho=1/N$, distinct integer prefix sums imply that the attractor has positive Lebesgue measure.

\begin{proof}[Proof of Theorem~\ref{thm:critical}]
Since the digits are commensurable, we may write
\[
\mathscr D=\alpha+\beta\mathscr E,
\qquad
\beta>0,\qquad
\mathscr E=\{0=e_0<e_1<\cdots<e_{N-1}\}\subset\Z,
\]
where $d_i=\alpha+\beta e_i$ and $\alpha=d_0$.
Put
$
\widetilde\mu=\mu_{1/N,\mathscr E},
$
$
\widetilde K=\supp\widetilde\mu
$
and
$
T(y)=\frac{\alpha}{N-1}+\beta y.
$
For $\widetilde S_i(y)=(y+e_i)/N$, a direct calculation gives
\[
S_i(T(y))
=\frac{\alpha}{N-1}+\frac{\beta}{N}(y+e_i)
=T(\widetilde S_i(y)).
\]
Thus $S_i\circ T=T\circ\widetilde S_i$.
By uniqueness of the self-similar probability measure, we obtain $\mu=T_{\#}\widetilde\mu$ and $K=T(\widetilde K).$

All six conditions are preserved by this affine change, as follows.
Frame-measure existence is preserved by
Lemma~\ref{lem:affine}.
For every word $w$, the conjugacy gives
$S_w=T\circ\widetilde S_w\circ T^{-1}$, so exact overlaps
are preserved as well.
Moreover,
\[
\sum_{j=0}^{n-1}N^j\mathscr D
=
\alpha\frac{N^n-1}{N-1}
+\beta\sum_{j=0}^{n-1}N^j\mathscr E,
\]
and hence $
\#\left(\sum_{j=0}^{n-1}N^j\mathscr D\right)
=
\#\left(\sum_{j=0}^{n-1}N^j\mathscr E\right).
$
Since $|K|=\beta|\widetilde K|$, positivity of Lebesgue
measure is preserved. When $|\widetilde K|>0$, a change
of variables also gives
\[
T_{\#}\left(\frac{\ind_{\widetilde K}}{|\widetilde K|}\,dy\right)
=
\frac{\ind_K}{|K|}\,dx.
\]
Finally, using the same symbol sequence, set
$\widetilde X=\sum_{j\ge1}N^{-j}e_{\omega_j}$.
Then $X=T(\widetilde X)$.
Since $T$ is a Borel bijection, it follows that
$
\sigma(X)=\sigma(\widetilde X)
$
and
$
H(\omega_1\mid X)=H(\omega_1\mid\widetilde X).
$ 
We may therefore assume throughout the remainder of
the proof that
$\mathscr D\subset\Z$ and $\min\mathscr D=0.$
Set
$
M=\max\mathscr D,$
$L_{\mathscr D}=\frac{M}{N-1}$ and
$\mathscr D_n=\mathscr D+N\mathscr D+\cdots+N^{n-1}\mathscr D.
$
Then $K\subset[0,L_{\mathscr D}]$.

Proposition~\ref{prop:exact} gives
\textup{(i)}$\Rightarrow$\textup{(ii)}.
To prove \textup{(ii)}$\Leftrightarrow$\textup{(iii)},
observe that, for $w=(w_1,\ldots,w_n)\in D_N^n$,
\begin{equation}\label{eq:S_w}
S_w(x)=N^{-n}(x+b_w),
\qquad
b_w=\sum_{j=1}^nN^{n-j}d_{w_j}.
\end{equation}
By the definition of $\mathscr D_n$, we have
$
\mathscr D_n
=
\{b_w:w\in D_N^n\}.
$
Thus $\mathscr D_n$ consists of the distinct values
obtained from the $N^n$ words of length $n$.
For any $u,v\in D_N^n$, the formulas in
\eqref{eq:S_w} give
\[
S_u=S_v\quad\Longleftrightarrow\quad b_u=b_v.
\]
If the level-$n$ maps are all distinct, the $N^n$
values $b_w$ are distinct, so $\#\mathscr D_n=N^n$.
Conversely, if two distinct words define the same map,
their values of $b_w$ coincide, and hence
$\#\mathscr D_n<N^n$.
Therefore there is no exact overlap at level $n$
if and only if $\#\mathscr D_n=N^n$.
Applying this equivalence for every $n\ge1$ proves
\textup{(ii)}$\Leftrightarrow$\textup{(iii)}.

Assume \textup{(iii)}.
Each $b\in\mathscr D_n$ then corresponds to exactly one
word of length $n$.
Iterating the self-similarity equation, we obtain,
for every bounded interval $I=[s,t]$ with $s\le t$,
\[
\mu(I)
=N^{-n}\sum_{b\in\mathscr D_n}\mu(N^nI-b).
\]
If $\mu(N^nI-b)>0$, then $N^nI-b$ meets
$K\subset[0,L_{\mathscr D}]$, which implies
$
b\in[N^ns-L_{\mathscr D},N^nt].
$ Hence
\[
\mu(I)
=N^{-n}
\sum_{\substack{b\in J_n}}
\mu(N^nI-b),
\]
where $J_n:=\mathscr D_n\cap
[N^ns-L_{\mathscr D},N^nt]$. Since $\mathscr D_n\subset\Z$ and an interval
of length $L$ contains at most $L+1$ distinct integers, we have $\# J_n\leq N^n|I|+L_{\mathscr D}+1$. 
Since each summand is at most $1$,
\[
\mu(I)
\le N^{-n}\bigl(N^n|I|+L_{\mathscr D}+1\bigr)
=|I|+(L_{\mathscr D}+1)N^{-n}.
\]
Letting $n\to\infty$ gives $\mu(I)\le|I|$.
Covering Borel sets by countably many intervals now gives $\mu\le dx$. Since $\mu(K)=1$, it follows that $|K|\ge1>0$. This proves \textup{(iii)}$\Rightarrow$\textup{(iv)}.

Assume \textup{(iv)}. Since $K=\bigcup_{i=0}^{N-1}S_iK$ and each $S_i$ has ratio $1/N$, we have
\[
|S_iK|=\frac{|K|}{N}\qquad\text{and}\qquad\sum_{i=0}^{N-1}|S_iK|=|K|.
\]
The first-level sets are therefore pairwise disjoint up to Lebesgue null sets. Define $\nu=|K|^{-1}\ind_K\,dx$.
For every Borel set $B$, change of variables gives
\[
\frac1N(S_i)_{\#}\nu(B)
=
\frac{|K\cap S_i^{-1}(B)|}{N|K|}
=
\frac{|B\cap S_iK|}{|K|}.
\]
Summing over $i$ and using the almost-everywhere disjointness of the first-level sets, we obtain
\[
\frac1N\sum_{i=0}^{N-1}(S_i)_{\#}\nu(B)
=
\frac{|B\cap K|}{|K|}
=\nu(B).
\]
Thus $\nu$ satisfies the same self-similarity equation
as $\mu$. Uniqueness gives $\mu=\nu$, proving
\textup{(iv)}$\Rightarrow$\textup{(v)}.

Under \textup{(v)}, let $f\in L^2(\mu)$ and set
$g=|K|^{-1}f\ind_K$.
Then $g\in L^2(\R)$ and $f\mu=g\,dx$, so
$\widehat{f\mu}=\widehat g$.
By Plancherel's theorem,
\begin{align*}
\int_{\R}|\widehat{f\mu}(\xi)|^2\,|K|\,d\xi
&=|K|\,\int_{\R}|\widehat{g}(\xi)|^2\,d\xi=|K|\,\int_{\R}|g(\xi)|^2\,d\xi\\
&=
|K|\int_{\R}
\left|\frac{f(x)\ind_K(x)}{|K|}\right|^2\,dx\\
&=
\frac1{|K|}\int_K|f(x)|^2\,dx\\
&=\|f\|_{L^2(\mu)}^2.
\end{align*}
Hence $|K|\,d\xi$ is a Parseval frame measure,
proving \textup{(v)}$\Rightarrow$\textup{(i)}.

The implication \textup{(i)}$\Rightarrow$\textup{(vi)} follows from
Theorem~\ref{thm:general-digits}, whereas
\textup{(vi)}$\Rightarrow$\textup{(ii)} follows from
Proposition~\ref{prop:exact}. This completes the equivalence.
\end{proof}

\section{Separated measures and rational powers}\label{sec:rational}

For separated consecutive-digit systems, the spatial point determines every prefix. The conditional prefix probabilities are therefore indicator functions, and $R_n=1$. Thus the argument of
Section~\ref{sec:overlap} gives no decay. Unique coding also allows us
to regard coordinate characters as functions in $L^2(\mu)$.
We use these characters as test functions and compare the poles of their mask quotients at different scales.

The analytic arguments in this section assume $0<\rho<1/N$,
and we write $q=\rho^{-1}>N$. The pole-cover definitions and arithmetic results require only $q>1$.
For an integer multiplier, the quotient method of \cite{FW}
uses a single target pole class. A noninteger rational multiplier may require several target classes. We establish the covering criterion needed to prove Theorem~\ref{thm:rational}.

\subsection{Finite pole covers and the analytic criterion}

Fix $N\ge2$ and $q>1$, put $\zeta=\ee^{2\pi i/N}$, and define
\begin{equation}\label{eq:quotient-intro}
 Q_r(t)=\frac{1-\ee^{-2\pi it}}{1-\zeta^r\ee^{-2\pi it}},
 \qquad P_r=\frac rN+\Z,\qquad 1\le r<N.
\end{equation}
The real poles of $Q_r$ are the points of $P_r$, and they are simple. A \emph{finite pole cover} of $P_r$ is a finite family of pairs $(m_\ell,r_\ell)$, with $m_\ell\in\N$ and $1\le r_\ell<N$, such that
\begin{equation*}
 P_r\subset\bigcup_{\ell=1}^L q^{-m_\ell}P_{r_\ell}.
\end{equation*}
Each pair is called a \emph{channel}:$m_\ell$ specifies the change
of scale, and $r_\ell$ specifies the target pole class. Pairs with different exponents count as different channels.

We also use a pointwise comparison condition, denoted by
$\BD(\alpha,\sigma)$, where $\alpha,\sigma\ge0$. This condition
requires constants $R_0\ge1$ and $c_0,K_0>0$ such that, for every
$R\ge R_0$, there is a nonempty finite family
\[
 \mathcal F_R=\{(m_{R,\ell},r_{R,\ell}):1\le\ell\le K(R)\},
 \qquad K(R)\le K_0R^\alpha,
\]
with $m_{R,\ell}\in\N$ and $1\le r_{R,\ell}<N$, for which
\begin{equation}\label{eq:BD-intro}
 |Q_r(t)|>R\quad\Longrightarrow\quad
 \max_{1\le\ell\le K(R)}|Q_{r_{R,\ell}}(q^{m_{R,\ell}}t)|
 \ge c_0R^{-\sigma}|Q_r(t)|.
\end{equation}
This is required on the common domain
\[
 \R\setminus\left(P_r\cup\bigcup_{\ell=1}^{K(R)}q^{-m_{R,\ell}}P_{r_{R,\ell}}\right).
\]
The family may depend on $R$, but not on $t$. The comparison must
hold at every point of the common domain where $|Q_r(t)|>R$.

\begin{theorem}\label{thm:equivalence}
Let $N\ge2$, $q>1$, and $1\le r<N$. The following conditions are equivalent.
\begin{enumerate}
\item The class $P_r$ has a finite pole cover.
\item The condition $\BD(\alpha,\sigma)$ holds for some $\alpha,\sigma\ge0$.
\item The number $q$ has a rational power. Writing
\begin{equation}\label{eq:minimal-power}
 d=\min\{m\in\N:q^m\in\Q\},\qquad q^d=\frac ab,
 \qquad a,b\in\N,\quad\gcd(a,b)=1,
\end{equation}
we have
\begin{equation}\label{eq:source-condition}
 b\mid N,\qquad b\mid r,\qquad N\nmid ar.
\end{equation}
\end{enumerate}
Whenever these conditions hold, $\BD(0,0)$ holds and the minimum number of channels is $b$. A minimum cover is
\begin{equation}\label{eq:optimal-intro}
 \{(d,r_h):0\le h<b\},\qquad
 r_h\equiv a\left(\frac Nb h+\frac rb\right)\pmod N,
 \quad1\le r_h<N.
\end{equation}
If $b>1$, every minimum cover has only exponent $d$, and its target residues are those in \eqref{eq:optimal-intro}. Consequently, there is a source residue $1\le r<N$ with a finite pole cover if and only if the least rational power exists and satisfies \eqref{eq:arithmetic-global}.
\end{theorem}

Thus allowing the number of channels and the comparison constant to
depend on the threshold as in \eqref{eq:BD-intro} gives no additional
arithmetic cases. The theorem determines exactly when the finite-cover
method applies. Each of its equivalent conditions yields the following
frame-measure obstruction.

\begin{theorem}\label{thm:equiv-obstruction}
Let $N\ge2$ and $q>N$. If a source residue $1\le r<N$ satisfies
any of the conditions of Theorem~\ref{thm:equivalence}, then
$\mu=\mu_{1/q,N}$ admits no frame measure. More precisely, for every
Bessel measure $\Gamma$ for $\mu$, there is $\tau\in\R$ such that
\begin{equation}\label{eq:walsh-decay-main}
 \int |\widehat{e_\tau W_{s,r}\mu}|^2\,d\Gamma\longrightarrow0
 \quad(s\to\infty),
\end{equation}
where $W_{s,r}(\sum_{j\ge1}q^{-j}\omega_j)=\zeta^{r\omega_s}$ and
$\|e_\tau W_{s,r}\|_{L^2(\mu)}=1$.
\end{theorem}

We prove Theorems~\ref{thm:equivalence} and~\ref{thm:equiv-obstruction}
after establishing the following lemmas. We first extend the two-coordinate
estimate of \cite{FW} to a finite family of target classes.
We then determine exactly when such a family exists.

\subsection{Cyclic characters and modulation}
Write $\mu=\mu_{\rho,N}$ and $K=\supp\mu$, and recall
$\zeta=\ee^{2\pi i/N}$. For $r=0,\ldots,N-1$, define
\[
m_r(t)=\frac1N\sum_{k=0}^{N-1}\zeta^{rk}\ee^{-2\pi ikt},
\qquad t\in\R.
\]
Since $0<\rho<1/N$, Lemma~\ref{lem:coding} shows that
$\pi:\Omega\to K$ is a homeomorphism. Thus each $x\in K$
has a unique coding $x=\pi(\omega)$. For $s\in\N$ and
$r\in\{1,\ldots,N-1\}$, define the coordinate character
$W_{s,r}:K\to\C$ by
\[
W_{s,r}(\pi(\omega))=\zeta^{r\omega_s},
\qquad \omega=(\omega_1,\omega_2,\ldots)\in\Omega.
\]
For a finite set of coordinate indices $F\subset\N$, choose
$r_j\in\{1,\ldots,N-1\}$ for each $j\in F$ and define
\begin{equation}\label{eq:Walsh}
W_F=\prod_{j\in F}W_{j,r_j},
\qquad
W_F(\pi(\omega))=\prod_{j\in F}\zeta^{r_j\omega_j}.
\end{equation}
We omit the dependence on the chosen residues $(r_j)_{j\in F}$
from the notation $W_F$ and set $W_\varnothing=1$.

The product on the right of \eqref{eq:Walsh} depends on only
finitely many coordinates and is continuous on $\Omega$.
Since $\pi^{-1}$ is continuous, $W_F$ is continuous on $K$.
Moreover, $|W_F(x)|=1$ for every $x\in K$, and hence
$\|W_F\|_{L^2(\mu)}^2=\int_K |W_F(x)|^2\,d\mu(x)=1$.
We next give product formulas for $\widehat\mu$ and
$\widehat{W_F\mu}$ and describe the poles of $Q_r$. 

\begin{lemma}\label{lem:products}
Let $\mu=\mu_{\rho,N}$ with $0<\rho<1/N$. For a finite set
$F\subset\N$, choose $r_j\in\{1,\ldots,N-1\}$ for each $j\in F$,
and let $W_F$ be defined by \eqref{eq:Walsh}. Then, for every
$\lambda\in\R$,
\begin{equation}\label{eq:Walsh-product}
\widehat\mu(\lambda)
=
\prod_{j=1}^{\infty}m_0(\rho^j\lambda)
\qquad\text{and}\qquad
\widehat{W_F\mu}(\lambda)
=
\prod_{j\in F}m_{r_j}(\rho^j\lambda)
\prod_{j\in\N\setminus F}m_0(\rho^j\lambda).
\end{equation}
Moreover, the real zero set of $\widehat\mu$ is
\begin{equation}\label{eq:zero-set}
Z(\widehat\mu)
=
\bigcup_{j\ge1}\bigcup_{r=1}^{N-1}\rho^{-j}P_r.
\end{equation}

For each $r\in\{1,\ldots,N-1\}$, the quotient $m_r/m_0$,
defined where $m_0\ne0$, extends to the meromorphic function
$Q_r$ in \eqref{eq:quotient-intro}. This function has period $1$,
satisfies $Q_r(0)=0$, and has simple real poles precisely at
the points of $P_r$. Moreover, there are constants
$a_r,b_r,\delta_r>0$, independent of $p\in P_r$, such that
\begin{equation}\label{eq:pole-bounds}
\frac{a_r}{|h|}
\le |Q_r(p+h)|
\le \frac{b_r}{|h|}
\qquad
(p\in P_r,\ h\in\R,\ 0<|h|<\delta_r).
\end{equation}
\end{lemma}

\begin{proof}	
Fix an integer $k\ge\max(F\cup\{1\})$, and set $r_j=0$ for $j\notin F$. Independence gives
\[
 \E\left[\prod_{j=1}^k
 \zeta^{r_j\omega_j}\ee^{-2\pi i\lambda\rho^j\omega_j}\right]
 =\prod_{j=1}^k\E\left[
 \zeta^{r_j\omega_j}\ee^{-2\pi i\lambda\rho^j\omega_j}\right]
 =\prod_{j=1}^k m_{r_j}(\rho^j\lambda).
\]
The partial coding sums converge uniformly, and the integrands have modulus one. Letting $k\to\infty$ and applying dominated convergence proves the second identity in \eqref{eq:Walsh-product}; the first follows by taking $F=\varnothing$. Writing $z=\ee^{-2\pi it}$, the geometric sums give
\[
m_0(t)=\frac{1-z^N}{N(1-z)},\qquad
m_r(t)=\frac{1-z^N}{N(1-\zeta^rz)},
\]
with removable values supplied by the finite sums. In particular,
$m_0(t)=1$ for $t\in\Z$, and
$Z(m_0)=\bigcup_{r=1}^{N-1}P_r$.
A zero factor in the product for $\widehat\mu$ gives a zero of
$\widehat\mu$. To prove the converse, fix $\lambda$ for which all
factors are nonzero. The defining sum gives
\[
 |m_0(t)-1|\le\frac1N\sum_{k=0}^{N-1}2\pi k|t|
 =\pi(N-1)|t|.
\]
For all sufficiently large $j$, the numbers $m_0(\rho^j\lambda)$ lie in $|z-1|<1/2$,
where the logarithm with $\log1=0$ satisfies $|\log z|\le2|z-1|$.
Consequently,
\[
 \sum_{j\ge j_0}|\log m_0(\rho^j\lambda)|
 \le2\pi(N-1)|\lambda|\sum_{j\ge j_0}\rho^j<\infty.
\]
The tail product is the exponential of a convergent series and is
nonzero. The preceding finite product is also nonzero. This proves
\eqref{eq:zero-set}.

Next, fix $r\in\{1,\ldots,N-1\}$. Dividing the geometric
sums and using continuity at removable singularities gives
\[
\frac{m_r(t)}{m_0(t)}
=
\frac{1-\ee^{-2\pi it}}
{1-\zeta^r\ee^{-2\pi it}}
=
Q_r(t)
\qquad\text{when }m_0(t)\ne0.
\]
The expression on the right defines a meromorphic function
on $\C$. It has period $1$ and satisfies $Q_r(0)=0$.
Its denominator vanishes at a real point $t$ precisely
when $t\in P_r$. For $p=k+r/N\in P_r$, we have
$
Q_r(p+h)
=
\frac{1-\zeta^{-r}\ee^{-2\pi ih}}
{1-\ee^{-2\pi ih}}.
$
Since $1-\ee^{-2\pi ih}=2\pi ih+O(h^2)$ and $1-\zeta^{-r}\ee^{-2\pi ih}=1-\zeta^{-r}+O(h),$
we obtain
\[
Q_r(p+h)
=
\frac{1-\zeta^{-r}}{2\pi ih}+O(1)
\qquad(h\to0).
\]
The coefficient of $h^{-1}$ is nonzero, so every point of
$P_r$ is a simple pole. Moreover, the expression for
$Q_r(p+h)$ is independent of $k$, so the remainder is
uniform over $p\in P_r$. Put $c_r=|1-\zeta^{-r}|/(2\pi)>0$. The preceding expansion
allows us to choose $0<\delta_r<1/2$ such that
\[
\frac{c_r}{2}
\le
|hQ_r(p+h)|
\le
\frac{3c_r}{2}
\qquad
(p\in P_r,\ h\in\R,\ 0<|h|<\delta_r).
\]
Taking $a_r=c_r/2$ and $b_r=3c_r/2$ proves
\eqref{eq:pole-bounds}.
\end{proof}

To use the mask quotients $Q_r$, we need to avoid the zeros of $\widehat\mu$. The following lemma shows that this can be done by a suitable modulation.

\begin{lemma}\label{lem:modulation}
Let $\nu$ be a compactly supported probability measure on $\R$, and let $\Gamma$ be a Bessel measure for $\nu$. Then $\Gamma$ is locally finite and there is a countable set $E_\Gamma$ such that
\begin{equation}\label{eq:generic}
 \Gamma\bigl(\tau+Z(\widehat\nu)\bigr)=0
 \qquad(\tau\notin E_\Gamma).
\end{equation}
\end{lemma}

\begin{proof}
Since $\widehat\nu$ is continuous and $\widehat\nu(0)=1$, we may choose $\delta>0$ such that $|\widehat\nu(t)|\ge1/2$ whenever $|t|<\delta$. Let $C$ be a Bessel bound for $\Gamma$. For every $x\in\R$, applying the Bessel inequality to $e_x(y)=\ee^{2\pi ixy}$ gives
\[
\frac14\Gamma((x-\delta,x+\delta))
\le
\int_{\R}|\widehat\nu(\lambda-x)|^2\,d\Gamma(\lambda)
\le C.
\]
Since every compact interval can be covered by finitely many such intervals, $\Gamma$ is locally finite. Let $\mathcal A_\Gamma=\{\lambda\in\R:\Gamma(\{\lambda\})>0\}$
be the set of atoms of $\Gamma$. By local finiteness, each bounded interval contains only finitely many atoms of mass at least $1/k$, for every $k\in\N$. Hence $\mathcal A_\Gamma$ is countable.

On the other hand, since $\nu$ has compact support, $\widehat\nu$ extends to an entire function on $\C$. As $\widehat\nu(0)=1$, its zeros are isolated, and hence the real zero set $Z(\widehat\nu)$ is countable. Now set $E_\Gamma=\mathcal A_\Gamma-Z(\widehat\nu)$.
Since both $\mathcal A_\Gamma$ and $Z(\widehat\nu)$ are countable, so is $E_\Gamma$. If $\tau\notin E_\Gamma$, then
$(\tau+Z(\widehat\nu))\cap\mathcal A_\Gamma=\varnothing.$
Thus every point of $\tau+Z(\widehat\nu)$ has zero $\Gamma$-measure. Since $\tau+Z(\widehat\nu)$ is countable, countable additivity gives
	\[
	\Gamma(\tau+Z(\widehat\nu))=0,
	\]
	which proves \eqref{eq:generic}.
\end{proof}

With the modulation chosen in Lemma~\ref{lem:modulation}, the Fourier
energy of each character product can be written in terms of the
mask quotients $Q_r$.

\begin{lemma}\label{lem:moments}
Let $\Gamma$ be a Bessel measure for $\mu$ with bound $C$. Choose $\tau$ as in Lemma~\ref{lem:modulation}, and define
\begin{equation*}
 d\eta(\lambda)=|\widehat\mu(\lambda-\tau)|^2\,d\Gamma(\lambda).
\end{equation*}
Then $\eta(\R)\le C$ and, for every character product $W_F$ defined in \eqref{eq:Walsh},
\begin{equation}\label{eq:moment}
 \int\left|\prod_{j\in F}Q_{r_j}(\rho^j(\lambda-\tau))\right|^2\,d\eta(\lambda)
 =\mathcal E_\Gamma(e_\tau W_F)\le C.
\end{equation}
If $\Gamma$ has a frame lower bound $A$, the same integral is at least $A$.
\end{lemma}

\begin{proof}
Since $\|e_\tau\|_{L^2(\mu)}=1$, the Bessel inequality gives
\[
\eta(\R)
=
\int_{\R}|\widehat\mu(\lambda-\tau)|^2\,d\Gamma(\lambda)
=
\int_{\R}|\widehat{e_\tau\mu}(\lambda)|^2\,d\Gamma(\lambda)
=
\mathcal E_\Gamma(e_\tau)
\le C.
\]
By Lemma~\ref{lem:modulation}, the set $B=\tau+Z(\widehat\mu)$ has zero $\Gamma$-measure. For $\lambda\notin B$, we have
$\widehat\mu(\lambda-\tau)\ne0$. Hence, by the first identity in \eqref{eq:Walsh-product}, every factor $m_0(\rho^j(\lambda-\tau))$ is nonzero. Since modulation translates the Fourier variable, using the two product formulas in \eqref{eq:Walsh-product} yields
\begin{align*}
\widehat{e_\tau W_F\mu}(\lambda)
&=
\widehat{W_F\mu}(\lambda-\tau)\\
&=
\widehat\mu(\lambda-\tau)
\prod_{j\in F}
\frac{m_{r_j}(\rho^j(\lambda-\tau))}
{m_0(\rho^j(\lambda-\tau))}\\
&=
\widehat\mu(\lambda-\tau)
\prod_{j\in F}Q_{r_j}(\rho^j(\lambda-\tau)),\qquad\lambda\notin B.
\end{align*}
This identity holds for $\Gamma$-almost every $\lambda$. For the purpose of these integrals, we may assign the value zero to $Q_r$ at its real poles, since these poles can only be encountered when $\lambda\in B$, and
$\Gamma(B)=0$.

Taking absolute values, squaring, and integrating
against $\Gamma$, we obtain
\[
\mathcal E_\Gamma(e_\tau W_F)
=
\int_{\R}
|\widehat{e_\tau W_F\mu}(\lambda)|^2\,d\Gamma(\lambda)
=
\int_{\R}
\left|\prod_{j\in F}
Q_{r_j}(\rho^j(\lambda-\tau))\right|^2
\,d\eta(\lambda).
\]
Finally, since $|e_\tau W_F|=1$ on $\supp\mu$, we have $\|e_\tau W_F\|_{L^2(\mu)}=1$.
The Bessel inequality therefore gives
$\mathcal E_\Gamma(e_\tau W_F)\le C.$
If $\Gamma$ has a frame lower bound $A$, the same
identity also gives $\mathcal E_\Gamma(e_\tau W_F)\ge A$. This proves \eqref{eq:moment}.
\end{proof}

In particular, for $s\ne t$ and nonzero residues $r,u$, we have
\begin{equation}\label{eq:two-moment}
\int|Q_r(\rho^s(\lambda-\tau))|^2\,d\eta\le C\quad\text{and}\quad
\int|Q_r(\rho^s(\lambda-\tau))Q_u(\rho^t(\lambda-\tau))|^2\,d\eta\le C.
\end{equation}
The condition $s\ne t$ in \eqref{eq:two-moment} is needed because the
two quotients must come from distinct coordinates. These estimates
hold for every Bessel measure $\Gamma$, including measures with atoms.

\subsection{A fixed-threshold estimate}
The following version of the two-coordinate argument in \cite{FW}
allows finitely many target functions. We fix the comparison
threshold before taking the limits in the proof.

\begin{lemma}\label{lem:tail}
Let $(Y,\eta)$ be a finite measure space, and let $v_s:Y\to\C$ be measurable functions such that $v_s\to0$ $\eta$-almost everywhere as $s\to\infty$. Suppose there are constants $R_*,c_*,C_*>0$, integers $L,s_0\ge1$, and measurable functions $w_{s,\ell}:Y\to\C$, $s\ge s_0$, $1\le\ell\le L$, such that, for every $s\ge s_0$,
\begin{equation}\label{eq:tail-hyp}
\max_{1\le\ell\le L}|w_{s,\ell}|\ge c_*|v_s|
\quad\text{on }\{|v_s|>R_*\},\qquad
\int|v_sw_{s,\ell}|^2\,d\eta\le C_*,
\end{equation}
where the first inequality holds $\eta$-almost everywhere. Then
\[
\int_Y|v_s|^2\,d\eta\longrightarrow0\qquad(s\to\infty).
\]
\end{lemma}

\begin{proof}
Fix $s\ge s_0$ and put $E_s=\{|v_s|>R_*\}$.
By the first inequality in \eqref{eq:tail-hyp},
\[
c_*^2|v_s|^4
\le
|v_s|^2\max_{1\le\ell\le L}|w_{s,\ell}|^2
\le
\sum_{\ell=1}^L|v_sw_{s,\ell}|^2
\]
$\eta$-almost everywhere on $E_s$. Integrating and using the second inequality in
\eqref{eq:tail-hyp}, we obtain
\begin{equation}\label{eq:fourth}
\int_{E_s}|v_s|^4\,d\eta
\le
c_*^{-2}\sum_{\ell=1}^L
\int_{E_s}|v_sw_{s,\ell}|^2\,d\eta
\le
Lc_*^{-2}C_*.
\end{equation}

For $T\ge R_*$, the set $\{|v_s|>T\}$ is contained in $E_s$,
and $|v_s|^2\le T^{-2}|v_s|^4$ on this set. Hence
\begin{align*}
\int_{\{|v_s|>T\}}|v_s|^2\,d\eta
&\le
T^{-2}\int_{\{|v_s|>T\}}|v_s|^4\,d\eta\\
&\le
T^{-2}\int_{E_s}|v_s|^4\,d\eta\\
&\le
Lc_*^{-2}C_*T^{-2},
\end{align*}
where the last inequality follows from \eqref{eq:fourth}. This bound is independent of $s$, hence
for every $T\ge R_*$,
\begin{equation}\label{eq:tail-bound}
\sup_{s\ge s_0}
\int_{\{|v_s|>T\}}|v_s|^2\,d\eta
\le Lc_*^{-2}C_*T^{-2}.
\end{equation}

Finally, fix $T\ge R_*$. Since $v_s\to0$ $\eta$-almost everywhere, the functions
$|v_s|^2\ind_{\{|v_s|\le T\}}$ tend to zero
$\eta$-almost everywhere and are bounded by $T^2$. As $\eta(Y)<\infty$, dominated convergence gives
$
\int_{\{|v_s|\le T\}}|v_s|^2\,d\eta\longrightarrow0.
$
Splitting the integral over $\{|v_s|\le T\}$ and
$\{|v_s|>T\}$ and applying \eqref{eq:tail-bound}, we obtain
\[
\limsup_{s\to\infty}\int_Y |v_s|^2\,d\eta
\le
Lc_*^{-2}C_*T^{-2}.
\]
Letting $T\to\infty$ completes the proof.
\end{proof}

We next show that a finite pole cover is equivalent to a comparison of the mask quotients above a fixed threshold. This comparison will allow us to apply Lemma~\ref{lem:tail}.

\begin{lemma}\label{lem:comparison}
Let $q>1$ and $r\in\{1,\ldots,N-1\}$. Fix a finite family of channels $\{(m_\ell,r_\ell):1\le\ell\le L\}$, where $L\ge1$, $m_\ell\in\N$, and $r_\ell\in\{1,\ldots,N-1\}$. Then
$P_r\subset\bigcup_{\ell=1}^L q^{-m_\ell}P_{r_\ell}$ if and only if there are constants $R_*,c_*>0$ such that
\begin{equation}\label{eq:pole-comparison}
|Q_r(t)|>R_*
\quad\Longrightarrow\quad
\max_{1\le\ell\le L}|Q_{r_\ell}(q^{m_\ell}t)|
\ge c_*|Q_r(t)|
\end{equation}
for every $t$ in the common domain
$
\mathcal D
=
\R\setminus
\left(
P_r\cup\bigcup_{\ell=1}^L q^{-m_\ell}P_{r_\ell}
\right).
$
\end{lemma}

\begin{proof}
Suppose first that the family covers $P_r$. Using the constants in \eqref{eq:pole-bounds}, choose
\[
0<\delta<
\min\left\{
\frac12,\delta_r,
\min_{1\le\ell\le L}q^{-m_\ell}\delta_{r_\ell}
\right\},
\qquad
c_*=
\min_{1\le\ell\le L}
\frac{a_{r_\ell}}{q^{m_\ell}b_r}>0.
\]
For each $p\in P_r$, choose an index $\ell$ such that $q^{m_\ell}p\in P_{r_\ell}$.
If $t=p+h\in\mathcal D$ and $0<|h|<\delta$, then
$
q^{m_\ell}t=q^{m_\ell}p+q^{m_\ell}h,
$ and
$
0<|q^{m_\ell}h|<\delta_{r_\ell}.
$
Applying \eqref{eq:pole-bounds}, we have
\begin{equation}\label{eq:Qr}
|Q_{r_\ell}(q^{m_\ell}t)|
\ge
\frac{a_{r_\ell}}{q^{m_\ell}|h|}
\ge
c_*\frac{b_r}{|h|}
\ge
c_*|Q_r(t)|.
\end{equation}
The index $\ell$ may depend on $p$, but $\delta$ and $c_*$ are independent of $p$.

To choose the threshold, consider the compact set
\[
K_\delta
=
\{t\in[0,1]:\operatorname{dist}(t,P_r)\ge\delta\}.
\]
The function $Q_r$ is continuous on $K_\delta$, so it is bounded there. Since $Q_r$ has period $1$, the same bound holds whenever $\operatorname{dist}(t,P_r)\ge\delta$.
Choose $R_*>0$ larger than this bound. If $t\in\mathcal D$ and $|Q_r(t)|>R_*$, then
$t=p+h$ for some $p\in P_r$ with $0<|h|<\delta$. The preceding estimate \eqref{eq:Qr} therefore proves	\eqref{eq:pole-comparison}.

Conversely, suppose that \eqref{eq:pole-comparison} holds, and assume that some $p\in P_r$ is not covered. Then
$q^{m_\ell}p\notin P_{r_\ell},\;
(1\le\ell\le L).$ Thus each function $t\mapsto Q_{r_\ell}(q^{m_\ell}t)$ is continuous in a neighborhood of $p$. Since the family is finite, we can choose $0<\varepsilon<1/2$ and $M<\infty$ such that the interval $[p-\varepsilon,p+\varepsilon]$ contains no target pole and
\[
\max_{1\le\ell\le L}|Q_{r_\ell}(q^{m_\ell}t)|
\le M
\qquad(|t-p|\le\varepsilon).
\]
Since distinct points of $P_r$ are at least one unit apart
and $\varepsilon<1/2$, $p$ is the only source pole in
$[p-\varepsilon,p+\varepsilon]$. This interval also contains
no target pole. Hence $t\in\mathcal D$ whenever
$0<|t-p|<\varepsilon$.

Since $p$ is a pole of $Q_r$, we have $|Q_r(t)|\to\infty$ as $t\to p$. We may therefore choose such a $t$ with
$|Q_r(t)|>\max\{R_*,M/c_*\}.$
By \eqref{eq:pole-comparison},
\[
M
\ge
\max_{1\le\ell\le L}|Q_{r_\ell}(q^{m_\ell}t)|
\ge
c_*|Q_r(t)|
>
M,
\]
a contradiction. Thus every point of $P_r$ is covered.
\end{proof}

\begin{proof}[Equivalence of \textup{(i)} and \textup{(ii)} in Theorem~\ref{thm:equivalence}]
If $\BD(\alpha,\sigma)$ holds, fix one $R_*\ge R_0$. The corresponding family is finite, and $c_0R_*^{-\sigma}>0$.
Lemma~\ref{lem:comparison} therefore shows that this family covers
every pole in $P_r$. Conversely, the same lemma gives \eqref{eq:BD-intro} with
$\alpha=\sigma=0$ from a finite cover. The same family works for
every sufficiently large threshold.
\end{proof}

\begin{proof}[Proof of Theorem~\ref{thm:equiv-obstruction} under \textup{(i)} or \textup{(ii)}] Write $\rho=q^{-1}$ and $\mu=\mu_{\rho,N}$. Let $\Gamma$ be any Bessel measure for $\mu$ with bound $C$. Choose $\tau$ and $\eta$ as in Lemma~\ref{lem:moments}, so that
\[
B=\tau+Z(\widehat\mu),\qquad
\Gamma(B)=0,\qquad
d\eta(\lambda)=|\widehat\mu(\lambda-\tau)|^2\,d\Gamma(\lambda),
\qquad
\eta(\R)\le C.
\]
In particular, $\eta$ is finite and $\eta(B)=0$. By the preceding equivalence and Lemma~\ref{lem:comparison}, either assumption gives a fixed finite family of channels $\{(m_\ell,r_\ell):1\le\ell\le L\}$ and constants $R_*,c_*>0$ satisfying \eqref{eq:pole-comparison}. Put
$s_0=1+\max_{1\le\ell\le L}m_\ell.$ For $s\ge s_0$ and $\lambda\notin B$, define
\[
v_s(\lambda)=Q_r(\rho^s(\lambda-\tau)),
\qquad
w_{s,\ell}(\lambda)
=
Q_{r_\ell}(\rho^{s-m_\ell}(\lambda-\tau)),
\quad 1\le\ell\le L.
\]
Both $v_s(\lambda)$ and $w_{s,\ell}(\lambda)$ are finite
for $\lambda\notin B$. Indeed, for $\lambda\notin B$, \eqref{eq:zero-set} gives
\[
\rho^j(\lambda-\tau)\notin P_u
\qquad(j\ge1,\ 1\le u<N),
\]
and all the indices $s$ and $s-m_\ell$ are positive. We set $v_s$ and $w_{s,\ell}$ equal to zero on $B$.

We now verify the hypotheses of Lemma~\ref{lem:tail}. For each fixed $\lambda\notin B$, we have $\rho^s(\lambda-\tau)\to0$. Since $Q_r$ is continuous at zero and $Q_r(0)=0$, it follows that $v_s\to0$ $\eta$-almost everywhere. Next, for $s\ge s_0$ and $\lambda\notin B$, set $t=\rho^s(\lambda-\tau)$. Since
$q^{m_\ell}t
=
\rho^{s-m_\ell}(\lambda-\tau),
$
the source and target quotients in \eqref{eq:pole-comparison} are precisely $v_s(\lambda)$ and $w_{s,\ell}(\lambda)$.
Thus
\[
\max_{1\le\ell\le L}|w_{s,\ell}|
\ge c_*|v_s|
\quad\text{on }\{|v_s|>R_*\}
\]
$\eta$-almost everywhere. Finally, $m_\ell\ge1$ implies that $s$ and $s-m_\ell$ are distinct positive integers. The two-coordinate estimate
\eqref{eq:two-moment} therefore gives
\[
\int_{\R}|v_sw_{s,\ell}|^2\,d\eta
\le C
\qquad(s\ge s_0,\ 1\le\ell\le L).
\]
All hypotheses of Lemma~\ref{lem:tail} are satisfied with $C_*=C$. Hence, using \eqref{eq:moment} with $F=\{s\}$, we obtain
\[
\mathcal E_\Gamma(e_\tau W_{s,r})
=
\int_{\R}|v_s|^2\,d\eta
\longrightarrow0.
\]
This proves \eqref{eq:walsh-decay-main} for every Bessel measure $\Gamma$.

If $\Gamma$ were a frame measure with lower bound $A>0$, then $|e_\tau W_{s,r}|=1$ on $\supp\mu$ would give
\[
\mathcal E_\Gamma(e_\tau W_{s,r})
\ge
A\|e_\tau W_{s,r}\|_{L^2(\mu)}^2
=
A
\qquad(s\ge s_0),
\]
contradicting the preceding limit. Thus $\mu$ admits no frame measure. The case of condition \textup{(iii)} follows from the arithmetic equivalence proved below.
\end{proof}

To apply $\BD(\alpha,\sigma)$ directly, fix
$R_*\ge R_0$ and use the family of channels supplied at this
threshold. The constants in Lemma~\ref{lem:tail} are
$
L=K(R_*)\le K_0R_*^\alpha,
$
$
c_*=c_0R_*^{-\sigma}
$
and
$
C_*=C.
$
These constants and the starting index $s_0$ are independent
of $s$ and of the truncation level $T$. Thus
\eqref{eq:tail-bound} gives
\[
\sup_{s\ge s_0}
\int_{\{|v_s|>T\}}|v_s|^2\,d\eta
\le
K_0Cc_0^{-2}R_*^{\alpha+2\sigma}T^{-2},
\qquad T\ge R_*.
\]
As in the proof of Lemma~\ref{lem:tail}, first letting
$s\to\infty$ and then $T\to\infty$ gives
$\int|v_s|^2\,d\eta\to0$ for all $\alpha,\sigma\ge0$.
Choosing the comparison threshold to be $T$ would give a bound proportional to $T^{\alpha+2\sigma-2}$. Keeping the two thresholds separate avoids a restriction on $\alpha+2\sigma$.

\subsection{Arithmetic classification and the minimum number of channels}
We now determine when a finite pole cover exists and find the minimum number of channels. We begin by describing the arithmetic structure of a single channel. Throughout this subsection, we assume only $q>1$. Fix $1\le r<N$ and write $p_k=k+r/N$ for $k\in\Z$.

\begin{lemma}\label{lem:rational-channels}
Let $(m,u)$ be a channel, where $m\in\N$ and $1\le u<N$. If it covers $p_{k_0}$ for some $k_0\in\Z$, then $q^m\in\Q$. In this case, let $d,a,b$ be as in \eqref{eq:minimal-power}. Then $m=jd$ for some $j\ge1$. Moreover, the source indices covered by this channel form exactly one congruence class modulo $b^j$, namely,
\[
\{k\in\Z:q^mp_k\in P_u\}
=
k_0+b^j\Z.
\]
\end{lemma}
\begin{proof}
Since $(m,u)$ covers $p_{k_0}=k_0+r/N$ and
$P_u=u/N+\Z$, there is an integer $v_0$ such that
\begin{equation}\label{eq:q^np_k}
q^m p_{k_0}=q^m\left(k_0+\frac rN\right)=v_0+\frac uN.
\end{equation}
As $1\le r<N$, we have $Nk_0+r\ne0$. Hence
\[
q^m=\frac{Nv_0+u}{Nk_0+r}\in\Q.
\]

Let $d,a,b$ be as in \eqref{eq:minimal-power}. Dividing
$m$ by $d$, write $m=jd+s,$ $ j\ge0,$ $ 0\le s<d.$ Then
$
q^s=\frac{q^m}{(q^d)^j}\in\Q.
$
If $s>0$, this contradicts the minimality of $d$. Thus $s=0$, so $m=jd$ with $j\ge1$, and $q^m=a^j/b^j$.

Finally, for any $k\in\Z$, it follows from \eqref{eq:q^np_k} and the fact $p_k=k+r/N$ that
\[
q^mp_k
=
q^mp_{k_0}+q^m(p_k-p_{k_0})
=
v_0+\frac uN+\frac{a^j}{b^j}(k-k_0).
\]
Since $P_u=u/N+\Z$ and $v_0\in\Z$, the point $p_k$
is covered by $(m,u)$ if and only if
$(a^j/b^j)(k-k_0)$ is an integer. Using
$\gcd(a^j,b^j)=1$, we obtain
\begin{equation}\label{eq:congruence-channel}
q^mp_k\in P_u
\quad\Longleftrightarrow\quad
\frac{a^j}{b^j}(k-k_0)\in\Z
\quad\Longleftrightarrow\quad
k\equiv k_0\pmod{b^j}.
\end{equation}
Therefore the source indices covered by this channel
are exactly the integers in $k_0+b^j\Z$.
\end{proof}

We next use Lemma~\ref{lem:rational-channels} to derive the
arithmetic conditions required for a finite pole cover of $P_r$.

\begin{lemma}\label{lem:necessary}
If $P_r$ has a finite pole cover, then the least rational power $q^d=a/b$ in \eqref{eq:minimal-power} exists and satisfies \eqref{eq:source-condition}; that is, $b\mid N,$ $b\mid r$ and $N\nmid ar.$
\end{lemma}

\begin{proof}
Choose a finite pole cover of $P_r$ and discard all empty channels. By Lemma~\ref{lem:rational-channels}, $q$ has a rational power, so $d,a,b$ in \eqref{eq:minimal-power} are defined. The remaining cover can be written as
\[
\mathcal F
=
\{(j_\ell d,u_\ell):1\le\ell\le L\},
\qquad
j_\ell\in\N,\quad 1\le u_\ell<N.
\]

We first prove that $b\mid r$ and $b\mid N$. For every $k\in\Z$, the point $p_k=k+r/N$ is covered by some channel $(j_\ell d,u_\ell)$. Thus
\[
\frac{a^{j_\ell}}{b^{j_\ell}}
\left(k+\frac rN\right)
=
v+\frac{u_\ell}{N}
\]
for some $v\in\Z$, and hence
\[
a^{j_\ell}(Nk+r)
=
b^{j_\ell}(Nv+u_\ell).
\]
Since $\gcd(a^{j_\ell},b^{j_\ell})=1$, we obtain $b^{j_\ell}\mid Nk+r$, and therefore $b\mid Nk+r$. As $k$ was arbitrary, this holds for every $k\in\Z$.
Taking $k=0$ gives $b\mid r$, while taking $k=1$ gives $b\mid N+r$. Subtracting yields $b\mid N$.

It remains to prove that $N\nmid ar$. Suppose, to the
contrary, that $N\mid ar$, and put
\[
g=\gcd(a,N),\qquad
D=\frac Ng,\qquad
J=\max_{1\le\ell\le L}j_\ell.
\]
Since $N=gD$, the assumption $N\mid ar$ gives
$D\mid(a/g)r$. Together with $\gcd(a/g,D)=1$, this implies
$D\mid r$.

We now choose a single source index that is missed by
every channel. Since $g\mid a$ and $\gcd(a,b)=1$, we have
$\gcd(g,b^J)=1$. Thus there is an integer $k_*$ satisfying
\[
gk_*\equiv-\frac rD\pmod{b^J}.
\]
When $b=1$, this congruence is satisfied by every integer.
For the chosen $k_*$, write
$
gk_*+\frac rD=b^Jv_*
$
with $v_*\in\Z$. It follows that
\[
Nk_*+r
=
D\left(gk_*+\frac rD\right)
=
Db^Jv_*.
\]
Consequently, for every $1\le j\le J$,
\begin{equation*}
q^{jd}p_{k_*}
=
\frac{a^j}{b^j}\frac{Db^Jv_*}{N}
=
\frac{a^j}{g}b^{J-j}v_*
\in\Z,
\end{equation*}
since $g\mid a^j$ and $J-j\ge0$.

In particular, taking $j=j_\ell$ shows that every channel
in $\mathcal F$ sends the same source pole $p_{k_*}$ to
an integer. However, $P_{u_\ell}=u_\ell/N+\Z$ contains
no integer because $1\le u_\ell<N$. Thus no channel in
$\mathcal F$ covers $p_{k_*}$, contradicting the covering
assumption. This proves $N\nmid ar$.
\end{proof}

We next show that the arithmetic conditions in Lemma~\ref{lem:necessary} are also sufficient for a finite pole
cover. The following lemma gives an explicit cover of $P_r$
with $b$ channels.
\begin{lemma}\label{lem:explicit}
Suppose that $q^m=a/b$ in lowest terms for some $m\in\N$, and that $b\mid N$, $b\mid r$, $N\nmid ar.$ For $0\le h<b$, let $r_h\in\{0,\ldots,N-1\}$ be defined by
\begin{equation}\label{eq:explicit-residues}
r_h\equiv
a\left(\frac Nb h+\frac rb\right)\pmod N.
\end{equation}
Then the residues $r_h$ are distinct and nonzero, and the $b$ channels $(m,r_h)$ cover $P_r$; that is,
$P_r\subset\bigcup_{h=0}^{b-1}q^{-m}P_{r_h}.$ The exponent $m$ need not be minimal.
\end{lemma}

\begin{proof}
Write $N=bM$ and $r=bc$. Then \eqref{eq:explicit-residues} becomes
\begin{equation}\label{eq:rh}
r_h\equiv a(Mh+c)\pmod N,
\qquad 0\le h<b.
\end{equation}
We first show that each $r_h$ is nonzero. Fix
$h\in\{0,\ldots,b-1\}$. If $r_h=0$, then
$N\mid a(Mh+c)$. Multiplying by $b$ gives
\[
N\mid ab(Mh+c)=a(Nh+r),
\]
and hence $N\mid ar$, contrary to the assumption.
Thus $1\le r_h<N$ for every $h$.

To prove that the residues are distinct, suppose that $r_h=r_{h'}$ for some $0\le h,h'<b$. Subtracting their
defining congruences gives
\[
N\mid aM(h-h').
\]
Since $N=bM$, we obtain $b\mid a(h-h')$.
By $\gcd(a,b)=1$, this implies $b\mid h-h'$.
As $|h-h'|<b$, we must have $h=h'$.

It remains to verify the covering property. Let $p\in P_r$.
Since $p-r/N$ is an integer, division with remainder gives
a unique representation
$p=h+bz+\frac rN$, $0\le h<b,$ $z\in\Z.$ Using $q^m=a/b$, $N=bM$, and $r=bc$, we obtain
\[
q^mp
=
\frac ab\left(h+bz+\frac rN\right)
=
az+\frac{a(Mh+c)}{N}.
\]
It follows from \eqref{eq:rh} that
\[
q^mp-\frac{r_h}{N}
=
az+\frac{a(Mh+c)-r_h}{N}
\in\Z.
\]
Hence $q^mp\in P_{r_h}$, so the channel $(m,r_h)$ covers
$p$. Since $p\in P_r$ was arbitrary,
$
P_r\subset\bigcup_{h=0}^{b-1}q^{-m}P_{r_h}.
$
This completes the proof.	
\end{proof}

Taking $m=d$ in Lemma~\ref{lem:explicit} gives a cover with
$b$ channels. We now show that this number is minimal and
describe all covers attaining it when $b>1$.
\begin{lemma}\label{lem:optimal}
Suppose that $P_r$ has a finite pole cover, and let $d,a,b$
be as in \eqref{eq:minimal-power}. Then the minimum number
of channels in such a cover is $b$. Moreover, if $b>1$,
every cover with exactly $b$ channels consists, up to
ordering, of the channels $(d,r_h)$, $0\le h<b$, where
the residues $r_h$ are given by \eqref{eq:optimal-intro}.
\end{lemma}
\begin{proof}
Choose an arbitrary finite pole cover of $P_r$ and discard
all empty channels. By Lemma~\ref{lem:rational-channels},
the remaining family can be written as
\[
\mathcal F
=
\{(j_\ell d,u_\ell):1\le\ell\le L\},
\qquad
j_\ell\in\N,\quad 1\le u_\ell<N.
\]
Put $J=\max_{1\le\ell\le L}j_\ell.$ For each $\ell$, choose an index $k_\ell\in\Z$ such that
the $\ell$th channel covers $p_{k_\ell}$.
By \eqref{eq:congruence-channel}, this channel covers
precisely the points $p_k$ whose indices satisfy
\[
k\equiv k_\ell\pmod{b^{j_\ell}}.
\]
To count all channels on a common finite set, consider the source indices $0,\ldots,b^J-1$. Since $b^{j_\ell}$ divides $b^J$, exactly
$b^{J-j_\ell}$ of these indices satisfy the preceding congruence. Every one of the $b^J$ indices is covered
by at least one channel. Therefore,
\begin{equation}\label{eq:count}
b^J
\le
\sum_{\ell=1}^L b^{J-j_\ell}
\le
Lb^{J-1},
\end{equation}
where the last inequality follows from $j_\ell\ge1$.
Dividing by $b^{J-1}$ gives $L\ge b$. Thus the original
cover also has at least $b$ channels.

By Lemma~\ref{lem:necessary}, we have
$b\mid N$, $b\mid r$, and $N\nmid ar$.
Applying Lemma~\ref{lem:explicit} with $m=d$ gives a
cover with $b$ channels. Hence the lower bound is attained.

Now suppose that $b>1$ and a cover has exactly $b$ channels.
None can be empty, since otherwise fewer than $b$ nonempty
channels would cover $P_r$, contrary to the lower bound.
Thus $L=b$, and equality holds throughout \eqref{eq:count}.
If $j_\ell>1$ for some $\ell$, then
$b^{J-j_\ell}<b^{J-1},$
so the last inequality in \eqref{eq:count} would be strict.
Consequently, $j_\ell=1$ for every $\ell$, and all channels
have exponent $d$.

It remains to identify their target residues. Let $r_h$,
$0\le h<b$, be the residues in \eqref{eq:optimal-intro}.
By Lemma~\ref{lem:explicit}, these residues are distinct
and nonzero. Moreover,
\[
q^d\left(h+\frac rN\right)
=
\frac{a}{N}\left(\frac Nb h+\frac rb\right)
\in P_{r_h}
\qquad(0\le h<b).
\]
Since the classes $P_u$, $1\le u<N$, are pairwise disjoint, a channel with exponent $d$ can cover the source pole $h+r/N$ only if its target residue is $r_h$.
Thus every channel $(d,r_h)$ must occur in the cover.
There are $b$ such channels and the cover has exactly $b$ channels, so these are all its channels.	
\end{proof}

\begin{proof}[Completion of the proof of
Theorem~\ref{thm:equivalence}]

For a fixed source residue $r$, the equivalence of
\textup{(i)} and \textup{(ii)}, together with the conclusion
$\BD(0,0)$, was proved above using
Lemma~\ref{lem:comparison}. Lemma~\ref{lem:necessary} gives
\textup{(i)}$\Rightarrow$\textup{(iii)}. Conversely, under
\textup{(iii)}, Lemma~\ref{lem:explicit}, applied with $m=d$,
gives a finite pole cover of $P_r$. Thus the three conditions
are equivalent. Lemma~\ref{lem:optimal} shows that the minimum number of
channels is $b$. If $b>1$, it also shows that every cover
with $b$ channels uses only exponent $d$ and has the target
residues in \eqref{eq:optimal-intro}.

For the last assertion, suppose that $P_{r_0}$ has a finite
pole cover for some $1\le r_0<N$. By
Lemma~\ref{lem:necessary}, the least rational power
$q^d=a/b$ exists and satisfies
\[
b\mid N,\qquad b\mid r_0,\qquad N\nmid ar_0.
\]
Write $N=bM$ and $r_0=bc$, where $M,c\in\N$.
If $M\mid a$, then $bM\mid abc=ar_0$, contradicting
$N\nmid ar_0$. Hence $M\nmid a$. Together with $b\mid N$,
this gives \eqref{eq:arithmetic-global}.

Conversely, suppose that the least rational power $q^d=a/b$
exists and satisfies \eqref{eq:arithmetic-global}.
Put $M=N/b$. Since $M\nmid a$, we have $M>1$.
Choose $r_*=b$. Then $1\le r_*<N$ and $b\mid r_*$.
Moreover,
\[
N\mid ar_*
\quad\Longleftrightarrow\quad
bM\mid ab
\quad\Longleftrightarrow\quad
M\mid a,
\]
so $N\nmid ar_*$. Thus \eqref{eq:source-condition} holds
with $r=r_*$. Applying Lemma~\ref{lem:explicit} with $m=d$
gives a finite pole cover of $P_{r_*}=P_b$.
This proves the last assertion.	
\end{proof}

By Theorem~\ref{thm:equivalence}, condition \textup{(iii)}
implies \textup{(i)}. Thus the proof of
Theorem~\ref{thm:equiv-obstruction} under \textup{(i)}
also establishes its conclusion under \textup{(iii)}.

\begin{proof}[Proof of Theorem~\ref{thm:rational}]

Let $q=\rho^{-1}>N$ and write $q^m=a/b$ as in
\eqref{eq:rational-power}, without requiring $m$ to be minimal. Set $r_*=b$.
By the same argument as in the last part of the proof of
Theorem~\ref{thm:equivalence}, \eqref{eq:arithmetic-global}
gives $1\le r_*<N$ and $N\nmid ar_*$.

Since $b\mid N$ and $b\mid r_*$, Lemma~\ref{lem:explicit},
applied to the given rational power $q^m=a/b$, provides
a finite pole cover of $P_{r_*}$ using exponent $m$.
Theorem~\ref{thm:equiv-obstruction}, applied under condition
\textup{(i)}, then shows that $\mu_{\rho,N}$ admits no
frame measure.
\end{proof}

\begin{proof}[Proof of Corollary~\ref{cor:classification}]
Here the digits are consecutive. At $\rho=1/N$, uniform probability on $[0,1]$ satisfies \eqref{eq:ifs}: the $N$ first-level intervals partition $[0,1]$ up to endpoints, and the factor $1/N$ compensates for the spatial scaling. Uniqueness gives $\mu_{1/N,N}=\ind_{[0,1]}\,dx$, which has spectrum $\Z$. Theorem~\ref{thm:overlap} excludes $1/N<\rho<1$.

Let $q\ge2$ be an integer. If $q<N$, apply Theorem~\ref{thm:overlap}. If $q>N$ and $N\nmid q$, apply Theorem~\ref{thm:rational} with $m=b=1$ and $a=q$. In the remaining cases $N\mid q$, and Theorem~\ref{thm:DHL} gives a spectrum and hence a frame measure.
\end{proof}

\section{Examples and remaining parameters}\label{sec:remaining}

We first give supercritical examples with and without exact overlaps,
and then illustrate the role of the digit set at the critical contraction.
Next, we apply the entropy obstruction to three commensurable digits.
We conclude with rational-power examples and the consecutive-digit
parameters not covered by our results.	

\subsection{Supercritical and critical examples}
\leavevmode
The next example has nonconsecutive digits and no exact overlap. Its attractor is not an interval, but Theorem~\ref{thm:overlap} still applies.

\begin{example}\label{ex:supercritical-gap}
Let $\mathscr D=\{0,1,5\}$ and $\rho=3/8$, and put $\mu=\mu_{\rho,\mathscr D}$ and $K=\supp\mu$. The convex hull of $K$ is $[0,3]$, and
\[
S_0([0,3])=[0,9/8],\qquad
S_1([0,3])=[3/8,3/2],\qquad
S_2([0,3])=[15/8,3].
\]	
Since $K=\bigcup_{i=0}^2S_i(K)$, the interval $(3/2,15/8)$ is disjoint
from $K$. As $3\rho=9/8>1$,
Theorem~\ref{thm:overlap} shows that $\mu$ admits no frame measure. Moreover,
\eqref{eq:entropy-main} gives		
\[
H(\omega_1\mid X)\ge\log(9/8)>0.
\]	
The functions $f_n$ in \eqref{eq:endpoint-main} have $L^2(\mu)$-norm one and vanish outside
$[0,(3/8)^n)$. They satisfy $\mathcal E_\Gamma(f_n)\to0$ for every Bessel
measure $\Gamma$, by Lemmas~\ref{lem:translation-energy} and
\ref{lem:entropy-decay}.		
	
There is no exact overlap. Otherwise, subtracting the prefix sums of two distinct words gives $\sum_{j=1}^n(3/8)^j\delta_j=0$, with
$\delta_j\in\{0,\pm1,\pm4,\pm5\}$ not all zero.
Let $r$ be the largest index with $\delta_r\ne0$.
Multiplication by $8^r$ and reduction modulo $8$ give
$3^r\delta_r\equiv0\pmod8$.
Since $\gcd(3,8)=1$, this forces $8\mid\delta_r$, which is impossible.
Thus this example is covered by the entropy obstruction but not by
the exact-overlap criterion in Proposition~\ref{prop:exact}.		
\end{example}

The next example is a singular measure in the supercritical range. It also has an exact overlap.

\begin{example}\label{ex:singular-supercritical}
Let $N=3$, $\mathscr D=\{0,1,2\}$ and $\rho=1/2$, and put
$\mu=\mu_{1/2,\mathscr D}$.
The attractor is $[0,2]$, and $S_0\circ S_2=S_1\circ S_0$ is an exact overlap.
We show that $\mu$ is singular with respect to Lebesgue measure.
Indeed, with		
\[
m_0(t)=\frac13(1+\ee^{-2\pi it}+\ee^{-4\pi it}),
\]
the product formula gives, for every nonnegative integer $n$,
\[
\widehat\mu(2^n)=\widehat\mu(1)
=\prod_{j\ge1}m_0(2^{-j})\ne0
\qquad(n\ge0).
\]
No dyadic number $2^{-j}$ belongs to the zero classes
$1/3+\Z$ or $2/3+\Z$, and
$\sum_{j\ge1}|m_0(2^{-j})-1|<\infty$.
Thus the product is nonzero, and the Riemann--Lebesgue lemma
excludes absolute continuity.	
	
To conclude singularity, write $\mu=\mu_{\rm ac}+\mu_{\rm s}$ for its
Lebesgue decomposition. Each similarity preserves both
absolute continuity and singularity. Uniqueness of the
decomposition therefore gives
$\mu_{\rm ac}=\frac13\sum_{i=0}^2(S_i)_{\#}\mu_{\rm ac}$.
If $\mu_{\rm ac}(\R)>0$, its normalization would be the
self-similar probability measure for the same system, forcing
$\mu$ to be absolutely continuous by uniqueness. Hence
$\mu_{\rm ac}=0$. Finally, $1/2>1/3$, so
Theorem~\ref{thm:overlap} shows that $\mu$ admits no frame measure.		
\end{example}

An integer digit set containing exactly one representative of each residue class modulo $N$ satisfies \eqref{eq:critical-cardinality}.
Indeed, reduction modulo $N$ determines the first digit in an integer
prefix sum. Subtracting that digit and dividing by $N$ determines
the next digit, and the argument continues by induction. Theorem~\ref{thm:critical} therefore gives a frame
measure at $\rho=1/N$.	

\begin{example}\label{ex:critical-three}
At $\rho=1/3$, compare the digit sets $\{0,1,5\}$ and $\{0,1,3\}$. The first contains one representative of each residue class modulo $3$. By Theorem~\ref{thm:critical}, its measure is normalized Lebesgue measure on its attractor $K$ and admits the Parseval frame measure $|K|\,d\lambda$. For the second set, the indexed maps satisfy
\[
S_0\circ S_2(x)=\frac{x}{9}+\frac13=S_1\circ S_0(x).
\]
Thus there is an exact overlap. Theorem~\ref{thm:critical} implies that
this measure has positive conditional entropy and admits no frame
measure. Hence two commensurable, nonconsecutive digit sets can give
different answers at the same critical contraction.		
\end{example}

For integer digits, it is enough to check
\eqref{eq:critical-cardinality} at one explicitly bounded level.
The following criterion applies even when the digits do not form
a complete residue system modulo $N$.

\begin{proposition}\label{prop:critical-carry}
	Let $\mathscr D\subset\Z$ consist of $N\ge2$ distinct digits
	with $\min\mathscr D=0$, and put
	$M=\max\mathscr D$ and
	$J=\left\lfloor M/(N-1)\right\rfloor$.
	Then
	$\mu_{1/N,\mathscr D}$ admits a frame measure if and only if
	$
	\#\bigl(\mathscr D+N\mathscr D+\cdots+N^{2J}\mathscr D\bigr)
	=N^{2J+1}.
	$	
\end{proposition}

\begin{proof}
An exact overlap at level $n$ is equivalent to a relation
$
\sum_{j=0}^{n-1}N^j\delta_j=0,
$ 
$
\delta_j\in\mathscr D-\mathscr D,
$
in which not all $\delta_j$ are zero.
We show that any such relation implies one of length at most $2J+1$.		
	
Choose such a relation with $n$ minimal.
Then $\delta_{n-1}\ne0$, since otherwise the last
coefficient could be omitted.
Set $s_0=0$ and, for $1\le k<n$, define
$
s_k
=N^{-k}\sum_{j=0}^{k-1}N^j\delta_j
=-\sum_{j=k}^{n-1}N^{j-k}\delta_j.
$
The second expression shows that $s_k$ is an integer. It records the
contribution of the first $k$ digit differences after division by $N^k$.
Since $|\delta_j|\le M$, the first expression gives
\[
|s_k|
\le N^{-k}M\sum_{j=0}^{k-1}N^j
=\frac{M}{N-1}(1-N^{-k})
<\frac{M}{N-1}.
\]
Thus $s_0,\ldots,s_{n-1}$ are integers in $[-J,J]$.

We claim that these integers are distinct.
If $s_r=s_t$ for some $0\le r<t<n$, then
\[
\sum_{j=0}^{r-1}N^j\delta_j
+\sum_{j=t}^{n-1}N^{j-t+r}\delta_j
=N^r(s_r-s_t)=0.
\]
This relation is obtained by deleting
$\delta_r,\ldots,\delta_{t-1}$ and has length
$n-t+r<n$. It is still nontrivial, since its last
coefficient is $\delta_{n-1}\ne0$.
This contradicts the minimality of $n$.
Hence $s_0,\ldots,s_{n-1}$ are distinct.
There are only $2J+1$ integers in $[-J,J]$, so
$n\le2J+1$.	
	
Since $0\in\mathscr D$, any equality between the digit sums of two distinct words of length $n\le2J+1$ extends to one of length
$2J+1$ by adding zero digits to both sums.
Consequently, there are no exact overlaps at any level
if and only if
$
\#\bigl(\mathscr D+N\mathscr D+\cdots+N^{2J}\mathscr D\bigr)
=N^{2J+1}.
$
The frame-measure characterization now follows from
Theorem~\ref{thm:critical}.		
\end{proof}

\begin{example}\label{ex:critical-product}
Let $\rho=1/4$, $\mathscr D=\{0,1,8,9\}$, and
$\mu=\mu_{\rho,\mathscr D}$.
These digits do not form a complete residue system modulo $4$.
Nevertheless, the distinct-prefix condition in
Proposition~\ref{prop:critical-carry} holds. Indeed, every digit in $\mathscr D$ has a unique representation
$a+8b$, where $a,b\in\{0,1\}$. Hence		
\[
\sum_{j=0}^{n-1}4^j(a_j+8b_j)
=
\sum_{j=0}^{n-1}2^{2j}a_j
+
\sum_{j=0}^{n-1}2^{2j+3}b_j.
\]
The exponents in the first sum are even, whereas those in
the second are odd. Thus all the powers of $2$ appearing
here are distinct. Uniqueness of finite binary expansions
gives
$
\#\bigl(\mathscr D+4\mathscr D+\cdots+4^{n-1}\mathscr D\bigr)
=4^n$, $n\ge1.$
Proposition~\ref{prop:critical-carry} therefore shows that
$\mu$ admits a frame measure.

In fact, the attractor is
$
K=[0,1]\cup[2,3].
$
Indeed, this set satisfies
$
\frac{K}{4}\cup\frac{K+1}{4}=[0,1]$ and
$\frac{K+8}{4}\cup\frac{K+9}{4}=[2,3],$
so the assertion follows from uniqueness of the compact
invariant set. Since $|K|=2$,
Theorem~\ref{thm:critical} gives
$
d\mu(x)=\tfrac12\ind_K(x)\,dx.
$
By Plancherel's theorem,
$
d\Gamma(\lambda)=2\,d\lambda
$
is a Parseval frame measure for $\mu$.
Thus a complete residue system is sufficient, but not necessary,
for a frame measure at the critical contraction.			
\end{example}

The next example shows that the absence of exact overlaps
need not imply positive Lebesgue measure when the digits
are not commensurable.

\begin{example}\label{ex:irrational-critical}
Let $\rho=1/3$ and $\mathscr D=\{0,1,\sqrt2\}$.
There is no exact overlap. To see this, suppose that
$
\sum_{j=0}^{n-1}3^j \xi_j
=
\sum_{j=0}^{n-1}3^j \eta_j,
$ $\xi_j,\eta_j\in\mathscr D.$
Write
$
\xi_j-\eta_j=a_j+\sqrt2\,b_j,
$ $a_j,b_j\in\{-1,0,1\}.
$
Since $\sqrt2$ is irrational, we obtain
$
\sum_{j=0}^{n-1}3^j a_j=0$ and $\sum_{j=0}^{n-1}3^j b_j=0.
$
Reduction modulo $3$ gives $a_0=b_0=0$.
Dividing both identities by $3$ and repeating the argument,
we obtain $a_j=b_j=0$ for every $j$.
Thus $\xi_j=\eta_j$ for every $j$, so the system has no exact
overlaps. Equivalently,
$
\#\bigl(\mathscr D+3\mathscr D+\cdots+3^{n-1}\mathscr D\bigr)
=3^n$, $n\ge1$.		
	
On the other hand, let $F\subset\R^2$ be the attractor
of the maps
\[
z\longmapsto\frac{z+v}{3},
\qquad v\in\{(0,0),(1,0),(0,1)\}.
\]
The attractor $K$ of the original system satisfies
$
K=\{x+\sqrt2\,y:(x,y)\in F\}.
$
Since the projection parameter $\sqrt2$ is irrational,
Kenyon's projection theorem~\cite{Kenyon1997} gives $|K|=0$.
Thus distinct finite digit sums do not guarantee positive
Lebesgue measure for arbitrary real digits.
This example does not determine whether
$\mu_{1/3,\mathscr D}$ admits a frame measure.		
\end{example}

\subsection{An application to three commensurable digits}\label{subsec:three-digits}

\leavevmode
Suppose that a frame measure exists. Theorem~\ref{thm:general-digits}
then gives $H(\omega_1\mid X)=0$, so the first symbol is determined
by $X$ almost surely. Iterating the tail map recovers every symbol,
as shown below. We can therefore use functions of individual
symbols as test functions on the attractor without assuming strong
separation. We begin with the case in which the mask is bounded
away from zero.

\begin{lemma}\label{lem:bounded-mask}
Let $\mathscr D=\{d_0<\cdots<d_{N-1}\}\subset\R$, where $N\ge2$, and let
$0<\rho<1$. Put
\[
M_0(t)=\frac1N\sum_{j=0}^{N-1}\ee^{-2\pi i d_jt}.
\]
If $\inf_{t\in\R}|M_0(t)|>0$, then $\mu_{\rho,\mathscr D}$ admits no frame
measure. Consequently, existence of a frame measure implies
$\inf_{t\in\R}|M_0(t)|=0$.	
\end{lemma}
\begin{proof}
Suppose that $\Gamma$ is a frame measure for $\mu=\mu_{\rho,\mathscr D}$
with bounds $A,C$. By Theorem~\ref{thm:general-digits},
$H(\omega_1\mid X)=0$. 
As in the proof of Lemma~\ref{lem:singular-components}, choose
pairwise disjoint Borel sets $A_i$ on which $p_i=1$ and whose union
has full $\mu$-measure.
Set $\iota(x)=i$ on $A_i$ and $\iota(x)=0$ elsewhere.
Then $\iota:\R\to D_N$ is Borel and
$\omega_1=\iota(X)$ almost surely.
Define the Borel tail map
$
\mathcal R(x)=\rho^{-1}x-d_{\iota(x)}.
$
The tail identity gives $X_1=\mathcal R(X)$ almost surely. Applying the
same identities after each shift, and removing a countable union of
null sets, we obtain
\begin{equation}\label{eq:measurable-coordinates-general}
\omega_n=\iota(\mathcal R^{n-1}(X))\qquad(n\ge1)
\quad\PP\text{-almost surely}.
\end{equation}
Thus every function of one symbol is a measurable function of $X$, up to a null set.

Let $\zeta=\ee^{2\pi i/N}$ and define the coordinate test functions $W_n$ and the corresponding mask quotient by
\[
W_n(X)=\zeta^{\omega_n},\qquad
M_1(t)=\frac1N\sum_{j=0}^{N-1}\zeta^j\ee^{-2\pi i d_jt},\qquad
\mathcal Q(t)=\frac{M_1(t)}{M_0(t)}.
\]
Then $\|W_n\|_{L^2(\mu)}=1$ and $\mathcal Q(0)=0$.
Independence of the coordinates, followed by dominated convergence,
gives
\[
\widehat{W_n\mu}(\lambda)
=M_1(\rho^n\lambda)\prod_{j\ne n}M_0(\rho^j\lambda)
=\mathcal Q(\rho^n\lambda)\widehat\mu(\lambda).
\]
The hypothesis makes $\mathcal Q$ bounded on $\R$, and
$\mathcal Q(\rho^n\lambda)\to0$ for every fixed $\lambda$.
The upper frame bound applied to $1$ shows that
$d\eta(\lambda)=|\widehat\mu(\lambda)|^2\,d\Gamma(\lambda)$ is finite.
Therefore
\[
\mathcal E_\Gamma(W_n)
=\int_{\R}|\mathcal Q(\rho^n\lambda)|^2\,d\eta(\lambda)
\longrightarrow0
\]
by dominated convergence. This contradicts
$\mathcal E_\Gamma(W_n)\ge A\|W_n\|_{L^2(\mu)}^2=A$.	
\end{proof}

After an affine change, any three commensurable digits have the form
$\{0,u,v\}$ with integers $0<u<v$ and $\gcd(u,v)=1$.
Indeed, digits $\alpha+\beta\{0,u,v\}$, where $\beta>0$, give the affine
image of $\mu_{\rho,\{0,u,v\}}$ under
$x\mapsto\beta x+\rho\alpha/(1-\rho)$.
Lemma~\ref{lem:affine} therefore reduces the frame-measure problem to
these normalized digits. Fu and Wen~\cite{FuWen2017} classified the
spectral measures for bounded sequences of three-element integer
digit sets $\{0,a_n,b_n\}$ with $\gcd(a_n,b_n)=1$ for every $n$; the fixed set considered here is a special case.
The next corollary classifies frame measures when the reciprocal
contraction is an integer and gives two further nonexistence statements.	

\begin{corollary}\label{cor:three-digits-frame}
Let $0<u<v$ be integers with $\gcd(u,v)=1$, and put
$\mathscr D=\{0,u,v\}$.
\begin{enumerate}
\item For an integer $q\ge2$, the measure $\mu_{1/q,\mathscr D}$ admits a
frame measure if and only if
\begin{equation}\label{eq:three-digit-condition}
3\mid q
\quad\text{and}\quad
\{u\bmod3,v\bmod3\}=\{1,2\}.
\end{equation}
In this case, the measure is spectral.
\item If the residue condition in \eqref{eq:three-digit-condition} fails,
then $\mu_{\rho,\mathscr D}$ admits no frame measure for any $0<\rho<1$.
\item Let $0<\rho\le1/3$. If, for some integer $m\ge1$,
\begin{equation}\label{eq:three-digit-rational-power}
\rho^{-m}=\frac ab,\qquad
a,b\in\N,\qquad \gcd(a,b)=1,\qquad
b\mid3,\qquad \frac3b\nmid a,
\end{equation}
then $\mu_{\rho,\mathscr D}$ admits no frame measure.	
\end{enumerate}	
\end{corollary}
\begin{proof}
When $q=2$, Theorem~\ref{thm:overlap} excludes a frame measure because
$1/2>1/3$. Thus it suffices to consider $q\ge3$ in~\textup{(i)}.
We first determine the real zeros of
$
M_0(t)=\frac13\bigl(1+\ee^{-2\pi iut}+\ee^{-2\pi ivt}\bigr).
$
Put $P(z)=1+z^u+z^v$ and $\zeta=\ee^{2\pi i/3}$. If $|z|=1$ and
$P(z)=0$, then $1,z^u,z^v$ are three unit complex numbers with sum zero.
They must be the three cube roots of unity, so
$\{z^u,z^v\}=\{\zeta,\zeta^2\}$. Hence $z^{3u}=z^{3v}=1$.
Since $\gcd(u,v)=1$, we have $z^3=1$ and $z\ne1$.
It follows that $M_0$ has a real zero if and only if the residue condition
in \eqref{eq:three-digit-condition} holds. If it fails, continuity and
periodicity give $\inf_{\R}|M_0|>0$.
Lemma~\ref{lem:bounded-mask} proves~\textup{(ii)}.
	
Assume that the residue condition holds. Then
$Z(M_0)=P_1\cup P_2$, where $ P_r=\frac r3+\Z\quad(r=1,2)$.
These zeros are simple. Indeed, for $z=\zeta$ or $\zeta^2$,
\[
\operatorname{Re}\bigl(zP'(z)\bigr)
=\operatorname{Re}(uz^u+vz^v)=-\frac{u+v}{2}<0.
\]
Define
$ M_r(t)=M_0(t-r/3)$ and $\mathcal Q_r(t)=\frac{M_r(t)}{M_0(t)}$, $r=1,2.$
The quotient $\mathcal Q_r$ has simple poles precisely at $P_r$;
at the other zero class of $M_0$, both numerator and denominator have
simple zeros, so their quotient has a finite limit. We define the
quotient there by this limit. Also, $\mathcal Q_r(0)=0$.	
	
We prove~\textup{(iii)}. If the residue condition fails, part~\textup{(ii)} already gives the result.
Otherwise, the preceding descriptions of $M_0$ and $\mathcal Q_r$ apply.
Since $b\mid3$, either $b=1$ or $b=3$.
The latter would give $3/b=1\mid a$, contrary to
\eqref{eq:three-digit-rational-power}. Thus $b=1$ and $3\nmid a$.
Put $q=\rho^{-1}$, so that $q^m=a$, and let $r'\in\{1,2\}$ satisfy
$r'\equiv a\pmod3$. Multiplication by $a$ sends $P_1$ into $P_{r'}$.
The simple-pole expansions give	
\[
\lim_{t\to1/3}
\frac{\mathcal Q_{r'}(at)}{\mathcal Q_1(t)}
=\frac{M_0'(1/3)}{aM_0'(r'/3)}\ne0.
\]
Both functions in this ratio are $1$-periodic, and $\mathcal Q_1$ is
bounded off any fixed neighborhood of $P_1$. Hence there are constants
$R_*,\kappa>0$ such that
\begin{equation}\label{eq:three-digit-propagation}
|\mathcal Q_1(t)|>R_*
\quad\Longrightarrow\quad
|\mathcal Q_{r'}(at)|\ge\kappa|\mathcal Q_1(t)|
\end{equation}
whenever both displayed quotients are finite.

Suppose that $\Gamma$ is a frame measure for
$\mu=\mu_{\rho,\mathscr D}$, with bounds $A,C$.
The construction of \eqref{eq:measurable-coordinates-general} uses only
the frame assumption and therefore applies here.
Hence the functions	
\[
W_{s,r}(X)=\ee^{2\pi i r d_{\omega_s}/3},\qquad s\ge1,
\quad r=1,2,
\]
are well defined in $L^2(\mu)$ and have norm one.
Choose $\tau$ by Lemma~\ref{lem:modulation}, so that
$\Gamma(\tau+Z(\widehat\mu))=0$, and put
$
d\eta(\lambda)=|\widehat\mu(\lambda-\tau)|^2\,d\Gamma(\lambda).
$
The Bessel inequality for $e_\tau$ gives $\eta(\R)\le C$.
For $s>m$, coordinate independence gives
\begin{align*}
\widehat{e_\tau W_{s,1}\mu}(\lambda)
&=\mathcal Q_1(\rho^s(\lambda-\tau))\widehat\mu(\lambda-\tau),\\
\widehat{e_\tau W_{s,1}W_{s-m,r'}\mu}(\lambda)
&=\mathcal Q_1(\rho^s(\lambda-\tau))
\mathcal Q_{r'}(\rho^{s-m}(\lambda-\tau))
\widehat\mu(\lambda-\tau)
\end{align*}
outside $\tau+Z(\widehat\mu)$. The product formula
$\widehat\mu(\lambda-\tau)=\prod_{j\ge1}M_0(\rho^j(\lambda-\tau))$
shows that every denominator factor is nonzero there.
Set
$ v_s(\lambda)=\mathcal Q_1(\rho^s(\lambda-\tau)),$ 
$
w_s(\lambda)=\mathcal Q_{r'}(\rho^{s-m}(\lambda-\tau)),
$
and set both functions equal to zero on the common null set.
Here $s$ and $s-m$ are distinct positive indices. Both character
functions above have modulus one and hence $L^2(\mu)$-norm one. The frame inequalities therefore give
\begin{equation}\label{eq:three-digit-integral-bounds}
A\le\int|v_s|^2\,d\eta,
\qquad \int|v_sw_s|^2\,d\eta\le C.
\end{equation}	
We have $v_s\to0$ almost everywhere because $\mathcal Q_1(0)=0$.
Moreover, $\rho^{s-m}=a\rho^s$, so
\eqref{eq:three-digit-propagation} gives
$|w_s|\ge\kappa|v_s|$ on $\{|v_s|>R_*\}$ almost everywhere.
Lemma~\ref{lem:tail}, with one target function, now yields
$\int|v_s|^2\,d\eta\to0$, contradicting
\eqref{eq:three-digit-integral-bounds}.
This proves~\textup{(iii)}. The two distinct coordinates provide the
same lower and upper integral bounds as in
Section~\ref{sec:rational} and \cite{FW}.		
	
For~\textup{(i)}, necessity of the residue condition follows from	
\textup{(ii)}, and necessity of $3\mid q$ follows from
\textup{(iii)} with $m=b=1$ and $a=q$.
Conversely, under \eqref{eq:three-digit-condition}, take
$\mathscr L=\{0,q/3,2q/3\}$. The matrix
$
\frac1{\sqrt3}
\bigl(\ee^{2\pi i d\ell/q}\bigr)_{d\in\mathscr D,\,\ell\in\mathscr L}
$
is the unitary $3\times3$ Fourier matrix with its rows possibly permuted.
Thus $(q,\mathscr D,\mathscr L)$ is a Hadamard triple, and
$\mu_{1/q,\mathscr D}$ is spectral by
\cite{DutkayHaussermannLai2019}, as also follows from \cite{FuWen2017}.		
\end{proof}

The arithmetic restrictions in
\eqref{eq:three-digit-rational-power} specialize
\eqref{eq:arithmetic-global} to $N=3$. They reduce to $b=1$ and
$3\nmid a$, so part~\textup{(iii)} concerns integer powers of
$\rho^{-1}$. Theorem~\ref{thm:equivalence} explains why this method
does not treat a noninteger rational power: when $N=3$, a finite
cover of either source pole class requires the least rational power
to have denominator one. The quotients above have the same pole
classes, so the same restriction applies. Parameters outside this
criterion require another argument.

\begin{example}\label{ex:three-opposite}
Fix $\rho=1/6$ and take
\[
\mathscr D_+=\{0,2,7\},\qquad \mathscr D_-=\{0,1,3\}.
\]
By Corollary~\ref{cor:three-digits-frame},
$\mu_{1/6,\mathscr D_+}$ is spectral, whereas	
$\mu_{1/6,\mathscr D_-}$ admits no frame measure.
Both systems satisfy strong separation. Their support convex hulls are
$[0,7/5]$ and $[0,3/5]$, whose lengths are smaller than the respective
minimum digit gaps $2$ and $1$. Thus the first-level images of each hull
are disjoint. The coding is injective, so both measures have
$H(\omega_1\mid X)=0$.
Thus $H(\omega_1\mid X)=0$ is not sufficient for the existence of a
frame measure, even for strongly separated systems with the same
contraction and the same number of digits.		
\end{example}

\subsection{Fractional multipliers and remaining parameters}
\leavevmode
For four consecutive digits, Theorem~\ref{thm:rational} applies to rational powers with denominator greater than one, as the following example shows.

\begin{example}\label{ex:fractional}
Let $N=4$ and $q=9/2$, so $\rho=2/9<1/4$.
Theorem~\ref{thm:rational} applies with $a=9$, $b=2$, and $m=1$.
For the source class $P_2=1/2+\Z$,
\[
q\left(k+\frac12\right)=\frac{18k+9}{4}
\in
\begin{cases}
P_1,&k\text{ even},\\
P_3,&k\text{ odd}.
\end{cases}
\]
The pairs $(1,1)$ and $(1,3)$ specify the common exponent $m=1$	
and the two target classes. They form a pole cover, and
Theorem~\ref{thm:equivalence} shows that two such pairs are necessary.
No positive power of $q$ is an integer, since the reduced denominator
of $q^j$ is $2^j$. Thus this parameter is not covered by the integer-power criterion
in Theorem~\ref{thm:FW}. The same reasoning applies to every
$q=a/2>4$ with $a$ odd. Theorem~\ref{thm:rational} also applies to
$\rho=\sqrt{2/33}<1/4$, since $\rho^{-2}=33/2$.		
\end{example}

A denominator $b>1$ in Theorem~\ref{thm:rational} is possible only
when $N$ is composite. Indeed, if $N$ is prime and $b>1$ divides $N$,
then $b=N$, so the condition $N/b\nmid a$ cannot hold.
This also explains the restriction to $b=1$ in
Corollary~\ref{cor:three-digits-frame}\textup{(iii)}.	

We now summarize the results for consecutive digits. For fixed $N$, set
\begin{align*}
\mathcal S_N&=\{1/q:q\in\N,\ N\mid q\},\\
\mathcal A_N&=\left\{\rho\in(0,1/N):
\begin{array}{l}
\rho^{-m}=a/b\text{ for some }m,a,b\in\N,\ \gcd(a,b)=1,\\[-2pt]
b\mid N,\quad N/b\nmid a
\end{array}\right\},\\
\mathcal U_N&=(0,1/N)\setminus(\mathcal S_N\cup\mathcal A_N).
\end{align*}
Here $\mathcal S_N$ is the spectral parameter set,
$\mathcal A_N$ is the set excluded by the rational-power criterion,
and $\mathcal U_N$ contains the remaining subcritical parameters.
Together with $(1/N,1)$, these sets form a disjoint partition of $(0,1)$.
The following table summarizes our results and
Theorem~\ref{thm:DHL}.	
\begin{center}
\small
\begin{tabularx}{\textwidth}{@{}p{0.31\textwidth}X@{}}
\toprule
Parameter range & Frame-measure conclusion\\
\midrule
$\rho\in\mathcal S_N$ & Spectral; counting measure on a spectrum is a frame measure. This includes $\rho=1/N$.\\[4pt]
$1/N<\rho<1$ & No frame measure, by Theorem~\ref{thm:overlap}.\\[4pt]
$\rho\in\mathcal A_N$ & No frame measure, by Theorem~\ref{thm:rational}.\\[4pt]
$\rho\in\mathcal U_N$ & Non-spectral by Theorem~\ref{thm:DHL}; frame-measure existence is not decided here.\\
\bottomrule		
\end{tabularx}
\end{center}

Suppose that $q=\rho^{-1}>N$ has a rational power, and let $d$
be the least positive integer for which $q^d=a/b\in\Q$, in lowest terms.
Theorem~\ref{thm:equivalence} shows that no source class has a finite
pole cover exactly when $b\nmid N$, or when $b\mid N$ and $N/b\mid a$.
If $q$ has no rational power, no source class has a finite cover either.
After excluding $\mathcal S_N$, these alternatives describe
$\mathcal U_N$.
For instance, our results do not decide whether
$\mu_{2/5,2}$ or $\mu_{1/\sqrt6,2}$ admits a frame measure.
In the first case, $b=N=2$ makes $N/b=1$, so the required nondivisibility fails;
in the second, the least rational power is $q^2=6$ and	
\[
6\left(k+\frac12\right)=6k+3\in\Z,
\]
so the only source pole class is sent to the integer class.
The parameters $\rho=\sqrt2-1$ and $\rho=1/\ee$ also lie in
$\mathcal U_2$. Indeed, a rational value of $(1+\sqrt2)^m$ would be
fixed by conjugation, and hence would equal $(1-\sqrt2)^m$, which
is impossible. No positive integer power of $\ee$ is rational
because $\ee$ is transcendental.	

In the separated consecutive-digit case, injectivity gives
$H(\omega_1\mid X)=0$ and $R_n=1$, as noted at the start of
Section~\ref{sec:rational}. Thus both the entropy obstruction and the finite pole-cover
criterion leave $\mathcal U_N$ undecided.

\begin{conjecture}\label{conj:full}
Let $N\ge2$ and $0<\rho<1$. The measure $\mu_{\rho,N}$ admits a frame
measure if and only if $\rho=1/q$ for an integer $q$ divisible by $N$.	
\end{conjecture}

Theorem~\ref{thm:overlap} settles the entire supercritical range; together with spectrality at $\rho=1/N$, this reduces the conjecture to $0<\rho<1/N$.
Theorem~\ref{thm:equivalence} determines exactly when the finite-family
comparison \eqref{eq:BD-intro} is available. To settle the conjecture
on $\mathcal U_N$, one needs a different comparison or a new choice
of test functions.


\begin{thebibliography}{99}
\bibitem{Christensen} O. Christensen, \emph{An Introduction to Frames and Riesz Bases}, second ed., Applied and Numerical Harmonic Analysis, Birkh\"auser/Springer, Cham, 2016.
\bibitem{CoverThomas} T. M. Cover and J. A. Thomas, \emph{Elements of Information Theory}, second ed., Wiley-Interscience, Hoboken, NJ, 2006.
\bibitem{DHLai} X.-R. Dai, X.-G. He and C.-K. Lai, Spectral property of Cantor measures with consecutive digits, \emph{Adv. Math.} \textbf{242} (2013), 187--208.
\bibitem{DHL} X.-R. Dai, X.-G. He and K.-S. Lau, On spectral $N$-Bernoulli measures, \emph{Adv. Math.} \textbf{259} (2014), 511--531.
\bibitem{PLT} J. de Dios Pont, L. Liehr and M. A. Taylor, Cantor measures with odd base do not admit Fourier frames, preprint, \href{https://arxiv.org/abs/2607.08656}{arXiv:2607.08656}, 2026.
\bibitem{DHW} D. E. Dutkay, D. Han and E. Weber, Continuous and discrete Fourier frames for fractal measures, \emph{Trans. Amer. Math. Soc.} \textbf{366} (2014), 1213--1235.
\bibitem{DutkayHaussermannLai2019} D. E. Dutkay, J. Haussermann and C.-K. Lai, Hadamard triples generate self-affine spectral measures, \emph{Trans. Amer. Math. Soc.} \textbf{371} (2019), 1439--1481.
\bibitem{DL} D. E. Dutkay and C.-K. Lai, Uniformity of measures with Fourier frames, \emph{Adv. Math.} \textbf{252} (2014), 684--707.
\bibitem{FH} D.-J. Feng and H. Hu, Dimension theory of iterated function systems, \emph{Comm. Pure Appl. Math.} \textbf{62} (2009), 1435--1500.
\bibitem{FSW} X.-Y. Fu, Z.-J. Song and W.-J. Wang, A Walsh--quotient obstruction for Fourier frames on odd reciprocal-power Bernoulli convolutions, preprint, \href{https://arxiv.org/abs/2607.10547}{arXiv:2607.10547}, 2026.
\bibitem{FW} X.-Y. Fu and W.-J. Wang, Nonexistence of frame measures for separated uniform consecutive-digit Bernoulli convolutions, preprint, \href{https://arxiv.org/abs/2608.19676}{arXiv:2608.19676}, 2026.
\bibitem{FuWen2017} Y.-S. Fu and Z.-X. Wen, Spectrality of infinite convolutions with three-element digit sets, \emph{Monatsh. Math.} \textbf{183} (2017), no.~3, 465--485.
\bibitem{Hutchinson} J. E. Hutchinson, Fractals and self-similarity, \emph{Indiana Univ. Math. J.} \textbf{30} (1981), 713--747.
\bibitem{JP} P. E. T. Jorgensen and S. Pedersen, Dense analytic subspaces in fractal $L^2$-spaces, \emph{J. Anal. Math.} \textbf{75} (1998), 185--228.
\bibitem{Kenyon1997} R. Kenyon, Projecting the one-dimensional Sierpinski gasket, \emph{Israel J. Math.} \textbf{97} (1997), 221--238.
\bibitem{LW02} I. \L{}aba and Y. Wang, On spectral Cantor measures, \emph{J. Funct. Anal.} \textbf{193} (2002), 409--420.
\bibitem{Lai11} C.-K. Lai, On Fourier frame of absolutely continuous measures, \emph{J. Funct. Anal.} \textbf{261} (2011), 2877--2889.
\bibitem{LaiWang} C.-K. Lai and Y. Wang, Non-spectral fractal measures with Fourier frames, \emph{J. Fractal Geom.} \textbf{4} (2017), 305--327.
\bibitem{Strichartz} R. S. Strichartz, Mock Fourier series and transforms associated with certain Cantor measures, \emph{J. Anal. Math.} \textbf{81} (2000), 209--238.
\bibitem{Walters} P. Walters, \emph{An Introduction to Ergodic Theory}, Graduate Texts in Mathematics, vol. 79, Springer-Verlag, New York, 1982.
\bibitem{Williams} D. Williams, \emph{Probability with Martingales}, Cambridge Mathematical Textbooks, Cambridge University Press, Cambridge, 1991.
\end{thebibliography}
\end{document}